\documentclass[12pt,reqno]{amsart}
\usepackage{amsmath, amsfonts, epsfig, amssymb, hyperref}
\usepackage{graphicx, graphics}
\usepackage[usenames]{color}
{\theoremstyle{plain}
\newtheorem{theorem}{Theorem}[section]
\newtheorem{conjecture}[theorem]{Conjecture}
\newtheorem{proposition}[theorem]{Proposition}
\newtheorem{lemma}[theorem]{Lemma}
\newtheorem{corollary}[theorem]{Corollary}
}
{\theoremstyle{definition}

\newtheorem{remark}[theorem]{Remark}
}

\newcommand{\ga}{\gamma}
\newcommand{\lt}{\lambda}
\newcommand{\mt}{\rho}
\newcommand{\de}{\delta}
\newcommand{\al}{\alpha}
\newcommand{\be}{\beta}
\newcommand{\at}{\sigma}
\newcommand{\bt}{\tau}
\newcommand{\supp}[2]{{\rm Supp}(#1, #2)}
\newcommand{\p}[1]{P(#1)}
\newcommand{\pd}[1]{D(#1)}
\newcommand{\addone}[1]{(#1+1)}
\newcommand{\addn}[2]{#1 \cup (#2)}

\newcommand{\domleq}{\preceq}
\newcommand{\domgeq}{\succeq}
\def\bleu{\textcolor{blue}}

\def\rouge{\textcolor{red}}

\def\jaune{\textcolor{yellow}}

\newcommand{\N}{\mathbb{N}}

\newcommand{\skewedhook}{skewed hook}
\newcommand{\row}[1]{\mathrm{row}(#1)}
\newcommand{\col}[1]{\mathrm{col}(#1)}
\newcommand{\partition}[1]{{#1}^\circlearrowleft}

\def\carre{\jaune{\linethickness{\unitlength}\line(1,0){1}}}

\begin{document}
\title[Differences of products of Schur functions]{Some positive differences of products\\ of Schur functions}

\author{Fran{\c{c}}ois Bergeron}
\address{D\'epartement de Math\'ematiques\\
Universit\'e du Qu\'ebec \`a Montr\'eal\\
Montr\'eal (Qu\'ebec) H3C 3P8\\
Canada}
\email{bergeron.francois@uqam.ca}
\thanks{F. Bergeron is supported in part by NSERC-Canada and
  FQRNT-Qu\'ebec.}

\author{Peter McNamara}
\address{Laboratoire de Combinatoire et d'Informatique Math\'ematique\\
Universit\'e du \linebreak Qu\'ebec \`a Montr\'eal\\
Montr\'eal (Qu\'ebec) H3C 3P8\\
Canada}
\email{mcnamara@lacim.uqam.ca}

\begin{abstract}
The product $s_\mu s_\nu$ of two Schur functions is one of the most famous
examples of a Schur-positive function, i.e. a symmetric function which, when written
as a linear combination of Schur functions, has all positive coefficients.
We ask when expressions of the form $s_\lambda s_\rho - s_\mu s_\nu$ are Schur-positive.
This general question seems to be a difficult one, but
a conjecture of Fomin, Fulton, Li and Poon says that it is the case at least when $\lambda$ and $\rho$ are obtained from $\mu$ and $\nu$ by
redistributing the parts of $\mu$ and $\nu$ in a specific, yet natural, way.  We show
that their conjecture is true
in several significant cases.  We also formulate a skew-shape extension
of their conjecture, and prove several results which serve as evidence in favor of 
this extension.
Finally, we take a more global view by studying two classes of
partially ordered sets suggested by these questions.
\end{abstract}

\maketitle

{ \parskip=0pt\footnotesize \tableofcontents}
\parskip=8pt  

\section{Introduction}
It is not hard to become convinced that Schur functions form the most important linear basis of the ring of symmetric functions. Not only do they play a fundamental role in the theory of symmetric polynomials, but they also are of deep significance in representation theory, algebraic geometry, as well as in many areas of mathematical physics. 
At the heart of the manifold reasons for this deep role, one finds the {\em Littlewood-Richardson coefficients} $c^\theta_{\mu \nu}$ that appear as structure constants for the multiplication of Schur functions:
\begin{equation}\label{multiplication}
   s_\mu s_\nu = \sum_{\theta} c^\theta_{\mu \nu} s_\theta .
\end{equation}
Recall that Schur functions $s_\lambda$ are naturally indexed by partitions $\lambda$. As is usual, we consider these as decreasing sequences of positive integers
   $$\lambda_1\geq \lambda_2\geq \ldots \geq \lambda_k>0,$$
of some length $k$. In Schur functions, as opposed to Schur ``polynomials,'' the variables (omitted in the notation above) are 
infinite in number: ${\bf x}=x_1,x_2,\ldots$. To make the presentation as self contained as possible, we will recall a combinatorial description of Schur functions in Section \ref{background}.   
The Littlewood-Richardson rule gives an interpretation for the $c^\theta_{\mu \nu}$'s of (\ref{multiplication}) as the number of semistandard 
Young tableaux satisfying certain conditions.  In particular, it follows that $c^\theta_{\mu \nu}$ is a non-negative
integer.  The product of two Schur functions is thus one of the most famous examples of 
a \emph{Schur-positive} function, i.e. 
a symmetric function which, when written
as a linear combination of Schur functions, has all positive coefficients.

In general, we plan to investigate the Schur-positivity of expressions of the form 
\begin{equation}\label{diff_pairs}
s_\tau s_\sigma - s_\mu s_\nu,
\end{equation}
with similar statements for skew Schur functions.
In view of  (\ref{multiplication}), the Schur-positivity of (\ref{diff_pairs}) clearly translates into a set of inequalities for the respective
Littlewood-Richardson coefficients:
$$
   c^\theta_{\mu \nu}\leq c^\theta_{\tau \sigma}.
$$
It is well known that these are difficult problems stated in such generality.
We will concentrate on special cases of the following form. Given a pair of partitions $(\mu, \nu)$, what sort of operations
can we apply to this pair to yield another pair $(\tau, \sigma)$ such that (\ref{diff_pairs}) is 
Schur-positive?
Two such interesting operations are considered by 
Fomin, Fulton, Li and Poon in \cite{FFLPpr}. However, it is still conjectural that (\ref{diff_pairs}) actually is Schur-positive in all of the instances they consider.  One of their operations,  called the $\ast$-operation, has been studied in \cite{BBR04pr}, where the pertinent conjecture is shown to hold for general families and ``asymptotically.''  We will be interested in
their other operation, which we will refer to as the $\sim$-operation (``tilde operation''). 
It is defined as follows. 
For a pair $(\mu, \nu)$ of partitions, let 
\begin{equation}\label{decroissant}
    \ga_1 \geq \ga_2 \geq \cdots \geq \ga_{2p}
 \end{equation}
be the decreasing rearrangement of the $\mu_i$ and
$\nu_j$'s. Then, we set
\begin{equation}\label{tilde_par}
   (\mu,\nu)^\sim:=(\lt,\mt),
\end{equation}
with
  $$   \lt = (\ga_1, \ga_3, \ldots , \ga_{2p-1}), \qquad {\rm and}\qquad 
      \mt=(\ga_2, \ga_4, \ldots , \ga_{2p}).$$
The following conjecture serves as the starting point for our investigations:

\begin{conjecture}\label{fflp27}
\textnormal{(Fomin, Fulton, Li, Poon)}
For all $\mu$ and $\nu$, if $(\lt, \mt)=(\mu,\nu)^\sim$, then 
    $$s_\lt s_\mt - s_\mu s_\nu$$
  is Schur-positive.
\end{conjecture}

In other words, suppose we have a list $\ga: \ga_1 \geq \ga_2 \geq \cdots \geq 
\ga_m $ of parts which we wish to distribute between two 
partitions $\lt$ and $\mt$. At the (negligible)  cost of adding a zero part, we may as well suppose that 
$m$ is even.  
We are interested in knowing which 
distribution results in the ``largest'' expression $s_\lt s_\mt$.
The conjecture says that the distribution that is like the usual dealing of cards is optimal.  

\begin{remark}
At this stage, it is natural to ask what happens when we distribute between more
than two partitions.  More specifically, after some investigation, it is tempting to make the following seemingly more general conjecture:

\begin{conjecture}\label{mtilde}
For $m \geq 2$, and a sequence of partitions $\mu^1, \ldots, \mu^m$, let 
   $$\ga:\ga_1 \geq \ga_2 \geq \cdots \geq \ga_{mp} \geq 0$$ 
be the decreasing rearrangement
of the $\mu^i_j$'s (with zeroes appended, to make the length of $\ga$
divisible by $m$).  Define $m$ new partitions $\widetilde{\mu}^1$, $\ldots\ $, $\widetilde{\mu}^m$ by
\[
\widetilde{\mu}^i := (\ga_i, \ga_{i+m}, \ldots, \ga_{i+(p-1)m})
\]
for $i=1, \ldots, m$.  
Then 
    $$s_{\widetilde{\mu}^1}s_{\widetilde{\mu}^2}\cdots s_{\widetilde{\mu}^m} - s_{\mu^1}s_{\mu^2}\cdots s_{\mu^m}$$
is Schur-positive.  
\end{conjecture}

However, Conjecture~\ref{fflp27} actually implies Conjecture~\ref{mtilde}.  This is not obvious, 
but one elegant way to show the implication is to use the ``repainting'' argument of
\cite[Proof of Prop. 2.9]{FFLPpr}.
Because of this, we will generally not make special mention of the $m$-partition case in what follows.  
However, several of our proofs for the 2-partition case, 
particularly those in Section \ref{skewshape}, 
work equally well for the $m$-partition case. 
We refer the interested reader to Remark~\ref{skewmtilde} for further details.  
\end{remark}

Our exposition is organized as follows.  After recalling some background in Section \ref{background}, we consider, in each of Sections  
\ref{multfree} and \ref{nontrivial}, 
special cases of Conjecture~\ref{fflp27}.  In Section \ref{skewshape}, we formulate a skew shape generalization of
Conjecture~\ref{fflp27} and show it to be true in some non-trivial special cases, as well
as giving other evidence in its favor.  
In Section \ref{posets}, we discuss
two classes of partially ordered sets (posets)
that arise naturally in our investigations.  Finally, in Section \ref{explodedjt}
we show how one of these posets leads us to consider the idea of an ``exploded'' Jacobi-Trudi
matrix.

Before beginning in earnest, let us
make one comment about numerical evidence.
Using software of A. Buch \cite{BucSoftware} and J. Stembridge \cite{SteSoftware}, we have
verified Conjecture~\ref{fflp27} for all $\mu$ and $\nu$ such that $|\mu| + |\nu| \leq 35$.

\section{Background and Notation}\label{background}
The usual notation and notions relating to partitions of integers and symmetric functions are recalled here. Notice that we are using here the French outlook (right side up!) for diagrams. For a partition $\mu$, we respectively  write $|\mu|$, $\ell(\mu)$ and $\mu'$ for the {\em sum of the parts}, {\em number of parts}
and {\em conjugate} of $\mu$.
In particular, $\mu'_i$ denotes the length of the $i^{\rm th}$ column of the {\em Young diagram} of
$\mu$. Recall  that, for a given $\mu$, this diagram, which we also denote by $\mu$,  is the set of $1$ by $1$ squares, in $\N\times \N$, with upper right corners
   $$\{ (i,j)\ |\ 1\leq i \leq \mu_j\}.$$
Thus the diagram of $6411$ is geometrically represented as
 $$\begin{picture}(6,4)(0,0)
\put(0,0){\line(1,0){6}}
\put(0,1){\line(1,0){6}}
\multiput(0,1)(1,0){7}{\line(0,-1){1}}
\put(0,2){\line(1,0){4}}
\multiput(0,2)(1,0){5}{\line(0,-1){1}}
\put(0,3){\line(1,0){1}}
\multiput(0,3)(1,0){2}{\line(0,-1){1}}
\put(0,4){\line(1,0){1}}
\multiput(0,4)(1,0){2}{\line(0,-1){1}}
\end{picture}\ .$$
When $\alpha$ is contained in $\mu$ as a diagram, written $\alpha \subseteq \mu$, 
we can consider the {\em skew shape}, usually denoted $\mu/\alpha$, whose cells are those in the set difference  $\mu\setminus \alpha$.
Let us denote by
$$
   \mu \cup \nu:=( \ga_1, \ga_2, \cdots )
$$
 the partition obtained by taking 
the decreasing rearrangement
of the $\mu_i$ and $\nu_j$'s,  just as in (\ref{decroissant}). For instance, $\mu\cup \nu=55444211$, if $\mu=5444$ and $\nu=5211$.
For any partition $\mu$ and $i \geq \ell(\mu)$, we consider $\mu_i$ be be zero.
Also, we set
$$
   \mu + \nu:=(\mu_1 + \nu_1, \ldots, \mu_\ell + \nu_\ell),
$$
where $\ell = \max \{\ell(\mu), \ell(\nu) \}$.
As is usual when $|\mu| = |\nu|$,  
we write $\mu \domleq \nu$ to denote that $\mu$ is less than or equal to $\nu$ in
{\em dominance order}. This means that we have all of the inequalities
\begin{eqnarray}
   \mu_1&\leq &\nu_1\nonumber\\
    \mu_1+\mu_2&\leq &\nu_1+\nu_2\nonumber\\
    &\vdots\label{def_dom}\\
      \mu_1+\ldots+\mu_i&\leq &\nu_1+\ldots+\nu_i\nonumber\\
          &\vdots\nonumber
\end{eqnarray}
Conjugation is an anti-isomorphism with respect to dominance order; i.e. $\mu \domleq \nu
\Leftrightarrow \mu' \domgeq \nu'$.
We use $a^k$ in the list of parts of a partition to denote a sequence of $k$ parts of the same size $a$.  Thus, a partition of the form $(j, 1^k)$ has one part of size $j$ and $k$ parts of size $1$. Such shapes are called {\em hooks}.  
A {\em semistandard Young tableau} $t$ of shape $\mu$ is an integer {\em filling} 
$$t:\mu\longrightarrow \{1,2, \ldots \},$$
of the cells of $\mu$, such that values are {\em strictly increasing} up the columns of $\mu$, and {\em weakly  increasing} along rows. Thus
$$t(i,j)\leq t(i+1,j),\qquad {\rm and}\qquad t(i,j)<t(i,j+1),$$
whenever these statements should make sense. For integers $a< b <c$, the following are 
semistandard tableau
of shape $21$:
\begin{equation}\label{tabl3}
\begin{picture}(2,2)(0,0)
\put(0,0){\line(1,0){2}}\put(0.2,0.2){$a$}
\put(0,1){\line(1,0){2}}\put(0.2,1.2){$b$}\put(1.2,0.2){$a$}
\multiput(0,1)(1,0){3}{\line(0,-1){1}}
\put(0,2){\line(1,0){1}}
\multiput(0,2)(1,0){2}{\line(0,-1){1}}
\end{picture} 
\qquad
\begin{picture}(2,2)(0,0)
\put(0,0){\line(1,0){2}}\put(0.2,0.2){$a$}
\put(0,1){\line(1,0){2}}\put(0.2,1.2){$b$}\put(1.2,0.2){$b$}
\multiput(0,1)(1,0){3}{\line(0,-1){1}}
\put(0,2){\line(1,0){1}}
\multiput(0,2)(1,0){2}{\line(0,-1){1}}
\end{picture} 
\qquad
\begin{picture}(2,2)(0,0)
\put(0,0){\line(1,0){2}}\put(0.2,0,2){$a$}
\put(0,1){\line(1,0){2}}\put(0.2,1.2){$b$}\put(1.2,0.2){$c$}
\multiput(0,1)(1,0){3}{\line(0,-1){1}}
\put(0,2){\line(1,0){1}}
\multiput(0,2)(1,0){2}{\line(0,-1){1}}
\end{picture} 
\qquad
\begin{picture}(2,2)(0,0)
\put(0,0){\line(1,0){2}}\put(0.2,0.2){$a$}
\put(0,1){\line(1,0){2}}\put(0.2,1.2){$c$}\put(1.2,0.2){$b$}
\multiput(0,1)(1,0){3}{\line(0,-1){1}}
\put(0,2){\line(1,0){1}}
\multiput(0,2)(1,0){2}{\line(0,-1){1}}
\end{picture}\ . \end{equation}
One naturally extends the notion of semistandard tableaux to skew shapes.
Recall that we say that the \emph{content} of a tableau (or skew tableau) is $\lambda$ if the tableau contains $\lambda_i$ copies of $i$, for all $i$. 

To make our presentation self-contained, we now recall the usual basic definitions regarding symmetric functions, with notation following \cite{Mac95}. To each tableau $t$ of shape $\mu$ (or skew shape $\mu/\alpha$), we associate the monomial
  $${\bf x}_t:=\prod_{c\in \mu} x_{t(c)} \, .$$
Then, the Schur symmetric\footnote{It is not evident from this definition, but they are truly symmetric.} function can be defined as
\begin{equation}
   s_\mu({\bf x}):=\sum_{t}  {\bf x}_t, 
\end{equation}
where the sum runs over the set of  semistandard tableaux of shape $\mu$. The skew Schur function $s_{\mu/\alpha}$ is likewise defined.
Considering  
(\ref{tabl3}), it is easy to see that
$$s_{21}= \sum_{a < b} (x_a^2 x_b + x_a x_b^2) + \sum_{a < b < c} 2 x_a x_b x_c \,.$$
The {\em complete homogeneous} symmetric function, indexed by an integer $n$, is defined to be 
\begin{equation}\label{chsf}
   h_n: = s_n = \sum_{i_1 \leq \cdots \leq i_n} x_{i_1}\cdots x_{i_n}\,.
\end{equation}
We further set
$$h_\mu:=h_{\mu_1}\cdots h_{\mu_k},$$
for a partition $\mu$.  In a similar way, the {\em elementary} symmetric function
indexed by an integer $n$, is obtained by replacing the inequalities in (\ref{chsf}) by strict inequalities: 
$$e_n: = s_{1^n} = \sum_{i_1 < \cdots < i_n} x_{i_1}\cdots x_{i_n}\,.$$
Just as before, we set
$$e_\mu:=e_{\mu_1}\cdots e_{\mu_k}$$ 
for a partition $\mu$.
We will denote by $\omega$ the
the well-known involution defined by 
$\omega(h_\mu)= e_\mu$ or, alternatively, by $\omega(s_\mu) = s_{\mu'}$.  As one would hope,
$\omega(s_{\mu/\al}) = s_{\mu'/\al'}$.

 
We finish this background section by recalling one of the (many) classical combinatorial descriptions of the Littlewood-Richardson coefficients. The {\em  reading word} of a (skew) tableau is obtained by reading the entries of the tableau starting with the bottom row, from right to left, and going up the rows.  For instance, $11221312$ and $11221213$ are the respective reading words of the skew semistandard tableaux
$$\begin{picture}(4,4)(0,0)
\put(2,0){\line(1,0){2}}
\put(1,1){\line(1,0){3}}
\multiput(2,1)(1,0){3}{\line(0,-1){1}}
\put(0,2){\line(1,0){4}}
\multiput(1,2)(1,0){4}{\line(0,-1){1}}
\put(0,3){\line(1,0){2}}
\multiput(0,3)(1,0){3}{\line(0,-1){1}}
\put(0,4){\line(1,0){1}}
\multiput(0,4)(1,0){2}{\line(0,-1){1}}
\put(0.2,3.2){$2$}
\put(0.2,2.2){$1$}\put(1.2,2.2){$3$}
                            \put(1.2,1.2){$1$}\put(2.2,1.2){$2$}\put(3.2,1.2){$2$}
                                                        \put(2.2,0.2){$1$}\put(3.2,0.2){$1$}
\end{picture}
\qquad\qquad
\begin{picture}(4,4)(0,0)
\put(2,0){\line(1,0){2}}
\put(1,1){\line(1,0){3}}
\multiput(2,1)(1,0){3}{\line(0,-1){1}}
\put(0,2){\line(1,0){4}}
\multiput(1,2)(1,0){4}{\line(0,-1){1}}
\put(0,3){\line(1,0){2}}
\multiput(0,3)(1,0){3}{\line(0,-1){1}}
\put(0,4){\line(1,0){1}}
\multiput(0,4)(1,0){2}{\line(0,-1){1}}
\put(0.2,3.2){$3$}
\put(0.2,2.2){$1$}\put(1.2,2.2){$2$}
                            \put(1.2,1.2){$1$}\put(2.2,1.2){$2$}\put(3.2,1.2){$2$}
                                                        \put(2.2,0.2){$1$}\put(3.2,0.2){$1$}
\end{picture}\ .$$
Observe that both of these tableaux are of shape $\theta/\mu=4421/21$, and content  $\nu=431$. A {\em lattice permutation} is a sequence of positive integers $a_1a_2\cdots a_n$ such that in any initial factor $a_1a_2\cdots a_j$, the number of $i$'s is at least as great as the number of $i+1$'s, for all $i$.   A proof of the following assertion can be found in \cite{Ful97, ec2}.

\begin{enumerate}
\item[] {\em
{\bf Littlewood-Richardson Rule.}\quad The Littlewood-Richardson coefficient $c_{\mu \nu}^\theta$ is equal to the number of semistandard tableaux of shape $\theta /\mu$ and content $\nu$ whose reading word is a lattice permutation. }
\end{enumerate}

\noindent When a semistandard tableaux of shape $\theta/\mu$ has a lattice permutation as its  reading word, we say that we have an  {\em LR-filling} of the shape $\theta/\mu$. We observe that  $c_{21,431}^{4421}=2$ since we have exhibited above two LR-fillings of $4421/21$ of content $431$, and these are easily seen to be the only possibilities.

As it turns out, we can also use the Littlewood-Richardson rule to expand skew Schur functions
in terms of Schur functions, since
\begin{equation}
        s_{\theta/\mu} = \sum_{\nu} c_{\mu \nu}^\theta s_\nu\ .
\end{equation}
Therefore, the Schur expansion of $s_{\theta/\mu}$ can be read off from the contents of 
all the LR-fillings of $\theta/\mu$.  
It follows readily from the definition of skew Schur functions that 
\begin{equation}\label{skew_product}
  s_{\theta/\mu} s_{\pi/\nu}=s_{(\theta/\mu)*(\pi/\nu)},
 \end{equation}
with the skew shape $(\theta/\mu) * (\pi/\nu)$ constructed as follows:
\[
\begin{picture}(7,6)(0,0)
\put(0,4){\line(0,1){2}}
\put(1,3){\line(0,1){1}}
\put(3,1){\line(0,1){5}}
\put(5,0){\line(0,1){1}}
\put(6,1){\line(0,1){2}}
\put(7,0){\line(0,1){1}}
\put(6,1){\line(1,0){1}}
\put(5,0){\line(1,0){2}}
\put(3,1){\line(1,0){2}}
\put(1,3){\line(1,0){5}}
\put(0,4){\line(1,0){1}}
\put(0,6){\line(1,0){3}}
\put(0.7, 4.7){$\pi/\nu$}
\put(3.7, 1.7){$\theta/\mu$}
\end{picture}\ .
\]
Therefore, the coefficient of $s_\lambda$ in the product $s_{\theta/\mu} s_{\pi/\nu}$ is 
equal to the number of LR-fillings of the shape $(\theta/\mu) * (\pi/\nu)$ with content $\lambda$.
Let us underline in passing that the identity (\ref{skew_product}) is often applied from
the right-hand side to the 
left-hand side.  More concretely, it states that if a skew shape consists of disjoint pieces, then its
associated skew Schur function is simply the product of the skew Schur functions associated to the 
pieces.

\section{Multiplicity-Free Products of Schur Functions}\label{multfree}

We are now going to underline how results from \cite{FFLPpr} 
and \cite{Ste01} combine
to imply that Conjecture~\ref{fflp27} is true for several 
infinite classes of partition pairs.
Our discussion will be made simpler if we introduce, for any pair of partitions
$\mu$ and $\nu$, the notion of the \emph{support} $\supp{\mu}{\nu}$ as being the set of all partitions $\theta$ such that $c^\theta_{\mu \nu} \neq 0$.
In other words, $\supp{\mu}{\nu}$ consists of those partitions $\theta$ for which $s_\theta$ appears
with non-zero coefficient in the expansion of $s_\mu s_\nu$.
In particular, all elements $\theta$ of $\supp{\mu}{\nu}$ satisfy $|\theta| = |\mu| + |\nu|$.
The following result appears as \cite[Corollary 2.6]{FFLPpr}:

\begin{proposition}[FFLP]
Suppose we have $(\lt, \mt)=(\mu,\nu)^\sim$.  Then 
   $$\supp{\mu}{\nu} \subseteq \supp{\lt}{\mt}.$$
\end{proposition}

This immediately implies that Conjecture~\ref{fflp27} follows for all pairs of partitions $(\mu, \nu)$ 
satisfying $c^\theta_{\mu\nu} \leq 1$ for all $\theta$.  But the set of such multiplicity-free 
pairs has been completely
characterized by Stembridge in \cite{Ste01}.  Before stating his result, we need some 
terminology.  
A partition $\mu$ with at most one part size is said to be a \emph{rectangle}.  
If the Young diagram of $\mu$ then has either $k$ rows, or $k$ columns, 
we say that $\mu$ is a \emph{$k$-line rectangle}.  A \emph{fat hook} is
a partition with exactly two part sizes, and if it is possible to obtain
a rectangle by deleting a single row or column from the fat hook $\mu$,
then we say that $\mu$ is a \emph{near-rectangle}.  For example,
\setlength{\unitlength}{2mm}
$$
\begin{picture}(5,3)(0,0)
\put(0,0){\line(1,0){5}}
\put(0,1){\line(1,0){5}}
\multiput(0,1)(1,0){6}{\line(0,-1){1}}
\put(0,2){\line(1,0){5}}
\multiput(0,2)(1,0){6}{\line(0,-1){1}}
\put(0,3){\line(1,0){3}}
\multiput(0,3)(1,0){4}{\line(0,-1){1}}
\end{picture}
 \qquad
\begin{picture}(5,3)(0,0)
\put(0,0){\line(1,0){5}}
\put(0,1){\line(1,0){5}}
\multiput(0,1)(1,0){6}{\line(0,-1){1}}
\put(0,2){\line(1,0){3}}
\multiput(0,2)(1,0){4}{\line(0,-1){1}}
\put(0,3){\line(1,0){3}}
\multiput(0,3)(1,0){4}{\line(0,-1){1}}
\end{picture}
\qquad
\begin{picture}(3,5)(0,0)
\put(0,0){\line(1,0){3}}
\put(0,1){\line(1,0){3}}
\multiput(0,1)(1,0){4}{\line(0,-1){1}}
\put(0,2){\line(1,0){3}}
\multiput(0,2)(1,0){4}{\line(0,-1){1}}
\put(0,3){\line(1,0){3}}
\multiput(0,3)(1,0){4}{\line(0,-1){1}}
\put(0,4){\line(1,0){2}}
\multiput(0,4)(1,0){3}{\line(0,-1){1}}
\put(0,5){\line(1,0){2}}
\multiput(0,5)(1,0){3}{\line(0,-1){1}}
\end{picture}
\qquad
\begin{picture}(3,5)(0,0)
\put(0,0){\line(1,0){3}}
\put(0,1){\line(1,0){3}}
\multiput(0,1)(1,0){4}{\line(0,-1){1}}
\put(0,2){\line(1,0){3}}
\multiput(0,2)(1,0){4}{\line(0,-1){1}}
\put(0,3){\line(1,0){3}}
\multiput(0,3)(1,0){4}{\line(0,-1){1}}
\put(0,4){\line(1,0){1}}
\multiput(0,4)(1,0){2}{\line(0,-1){1}}
\put(0,5){\line(1,0){1}}
\multiput(0,5)(1,0){2}{\line(0,-1){1}}
\end{picture}$$
\setlength{\unitlength}{4mm}
are all 
near-rectangles.  Then, as shown in \cite{Ste01}:

\begin{theorem}[Stembridge]\label{stemb}
The product $s_\mu s_\nu$ is multiplicity-free if and only if
\begin{itemize}\itemsep=4pt
\item[(i)] $\mu$ or $\nu$ is a one-line rectangle, or
\item[(ii)] $\mu$ is a two-line rectangle and $\nu$ is a fat hook (or vice versa), or
\item[(iii)] $\mu$ is a rectangle and $\nu$ is a near-rectangle (or vice versa), or
\item[(iv)] $\mu$ and $\nu$ are rectangles.
\end{itemize}
\end{theorem}
While of a similar flavor, pairs of hooks are not multiplicity free since, for example,
$c^{321}_{21,21} = 2$.  However, we have the following result:

\begin{proposition}
Conjecture~\ref{fflp27} holds when $\mu$ and $\nu$ are both hooks.  
\end{proposition}
\begin{proof}[\bf Proof.]
Suppose $\mu = (\mu_1, 1^r)$ and 
$\nu=(\nu_1, 1^s)$.
Because of Theorem \ref{stemb}(i) above, we can assume that $r$ and $s$ are non zero, and without loss of generality that $\mu_1 \geq \nu_1\geq 2$.    
Let $u = \lceil \frac{r+s}{2} \rceil$ and $v=\lfloor \frac{r+s}{2} \rfloor$.
From the definition, we have $(\lambda,\rho)=(\mu,\nu)^\sim$ with
      $$\lt = (\mu_1, 1^{u}), \qquad {\rm and}\qquad \mt = (\nu_1, 1^{v}).$$ 
Now, fix a partition $\theta$ such that $c^\theta_{\mu \nu} \neq 0$.  By the Littlewood-Richardson
rule, $c^\theta_{\mu\nu}$ is equal to the number of Littlewood-Richardson fillings (LR-fillings)
of $\theta/\mu$ of content $\nu$.  We wish to construct an injection $f$ from the set of
LR-fillings of $\theta/\mu$, of content $\nu$, to the set of LR-fillings of $\theta/\lt$, of content
$\mt$.  If $r=s$ or $r=s+1$, then $f$ can just be the identity map. 
We will assume that $r < s$, with the case $r > s$ being similar.  Consider an
LR-filling $t$ of $\theta/\mu$ of content $\nu$.  Observe that the first and second rows of $\theta/\mu$ are
the only ones that can have length greater than 1.  
Since $t$ is an LR-filling, we see that
the $s - v$ highest entries of $t$ must all be together at the top of the first column.
Delete these entries and move all other entries of the first column of $t$ up $s - v$ 
squares.  The result, denoted $f(t)$, is clearly still an LR-filling, and it is not 
difficult to see that $f$ is an injection.  Furthermore, $f(t)$ has shape $\theta/\lt$
and content $\mt$, as required.
\end{proof}

\begin{figure}\center
\begin{picture}(9,9)(0,0)\put(-1.5,4){$t:$}
\put(0,3.5){\carre}
\put(0.2,3.2){1}
\put(0,4.5){\carre}
\put(0.2,4.2){3}
\put(0,5.5){\carre}
\put(0.2,5.2){4}
\put(0,6.5){\carre}
\put(0.2,6.2){5}
\put(0,7.5){\carre}
\put(0.2,7.2){6}
\put(1,1.5){\carre}
\put(1.2,1.2){2}
\put(4,0.5){\carre}
\put(4.2,0.2){1}
\put(5,0.5){\carre}
\put(5.2,0.2){1}
\put(0,0){\line(1,0){6}}
\put(0,1){\line(1,0){6}}
\multiput(0,1)(1,0){7}{\line(0,-1){1}}
\put(0,2){\line(1,0){2}}
\multiput(0,2)(1,0){3}{\line(0,-1){1}}
\put(0,3){\line(1,0){1}}
\multiput(0,3)(1,0){2}{\line(0,-1){1}}
\put(0,4){\line(1,0){1}}
\multiput(0,4)(1,0){2}{\line(0,-1){1}}
\put(0,5){\line(1,0){1}}
\multiput(0,5)(1,0){2}{\line(0,-1){1}}
\put(0,6){\line(1,0){1}}
\multiput(0,6)(1,0){2}{\line(0,-1){1}}
\put(0,7){\line(1,0){1}}
\multiput(0,7)(1,0){2}{\line(0,-1){1}}
\put(0,8){\line(1,0){1}}
\multiput(0,8)(1,0){2}{\line(0,-1){1}}
\end{picture}
\put(-3,4){$\longmapsto$}\qquad\qquad
\begin{picture}(7,7)(0,0)\put(-3,4){$f(t):$}
\put(0,5.5){\carre}
\put(0.2,5.2){1}
\put(0,6.5){\carre}
\put(0.2,6.2){3}
\put(0,7.5){\carre}
\put(0.2,7.2){4}
\put(1,1.5){\carre}
\put(1.2,1.2){2}
\put(4,0.5){\carre}
\put(4.2,0.2){1}
\put(5,0.5){\carre}
\put(5.2,0.2){1}
\put(0,0){\line(1,0){6}}
\put(0,1){\line(1,0){6}}
\multiput(0,1)(1,0){7}{\line(0,-1){1}}
\put(0,2){\line(1,0){2}}
\multiput(0,2)(1,0){3}{\line(0,-1){1}}
\put(0,3){\line(1,0){1}}
\multiput(0,3)(1,0){2}{\line(0,-1){1}}
\put(0,4){\line(1,0){1}}
\multiput(0,4)(1,0){2}{\line(0,-1){1}}
\put(0,5){\line(1,0){1}}
\multiput(0,5)(1,0){2}{\line(0,-1){1}}
\put(0,6){\line(1,0){1}}
\multiput(0,6)(1,0){2}{\line(0,-1){1}}
\put(0,7){\line(1,0){1}}
\multiput(0,7)(1,0){2}{\line(0,-1){1}}
\put(0,8){\line(1,0){1}}
\multiput(0,8)(1,0){2}{\line(0,-1){1}}
\end{picture}
\caption{Injection from LR-fillings of $\theta/\mu$ to LR-fillings of $\theta/\lambda$.}
\label{hook_ex}
\end{figure}
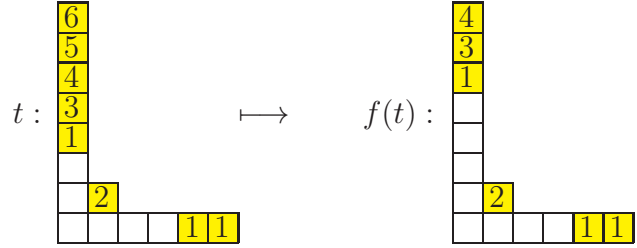

This argument is illustrated in Figure
\ref{hook_ex} for the shape $62111111$, which lies in both of the support sets $\supp{411}{311111}$
and  $\supp{411}{311111}^\sim=\supp{41111}{3111}$. We see here the resulting LR-filling $f(t)$, of content $\rho= 3111$, for a
given LR-filling $t$, of content $\nu=311111$.


\section{A $\ga$-Independent Special Case}\label{nontrivial}

For our next special case, referring to the notation of Conjecture~\ref{fflp27},
 we restrict to partition pairs of the following 
form:  for all $i$, $1 \leq i \leq p$,  we choose $\mu_i$ and $\nu_i$ so 
that $\{\mu_i, \nu_i\} = \{\ga_{2i-1},\ga_{2i} \}$ as multisets.
In other words, we take the parts of $\ga$ two at a time, giving one
to $\mu$ and the other to $\nu$.  Another way to express this is by saying that 
$\mu + \nu = \lt + \mt$.  

\begin{proposition}\label{specialcase}
Suppose that $(\mu, \nu)^\sim= (\lt, \mt)$, with $\mu + \nu = \lt + \mt$.
Then $$s_\lt s_\mt - s_\mu s_\nu$$ is Schur-positive.
\end{proposition}

\begin{remark}
In the setting of this proposition, 
we have $\mu + \nu = \lt + \mt$ and
$\mu \cup \nu = \lt \cup \mt$.  This is relevant to products of Schur functions because of
the following observation, which is not difficult to check: for any pair of partitions $\mu$, $\nu$,
let $\theta \in \supp{\mu}{\nu}$.  In dominance order, the following equation gives tight
bounds on $\theta$:
\[
\mu \cup \nu  \domleq \theta \domleq \mu+\nu.
\]
Therefore, the extreme elements of $\supp{\mu}{\nu}$ (in dominance order) are the same as those of
$\supp{\lambda}{\rho}$. 
\end{remark}

Our two main tools for proving Proposition \ref{specialcase}
will be the Jacobi-Trudi identity (\ref{jacobi_trudi})  and the Pl\"ucker
relations (\ref{plucker}), which we now state.  
In the form we need, the {\em Jacobi-Trudi identity} says that,
for any partition $\mu = (\mu_1, \mu_2, \ldots, \mu_p)$, we have
\begin{equation}\label{jacobi_trudi}
  s_\mu = \det(h_{\mu_i -i+j})^p_{i,j=1} .
\end{equation}
Here we allow $\mu_p=0$; and we set $h_0:=1$, and $h_j:=0$ if $j<0$. In other words, the matrix involved in (\ref{jacobi_trudi}) has $h_k$'s in the main diagonal, with indices equal to the parts of $\mu$:
\begin{equation}\label{JT_matrix}
    \begin{bmatrix}
                       h_{\mu_1} &\ast  &\ldots & \ast \\
                      \ast &h_{\mu_2}   &\ldots & \ast \\
                      \vdots &\vdots &\ddots &\vdots\\
                      \ast &\ast &\ldots & h_{\mu_k}
           \end{bmatrix} , 
\end{equation}
while rows are filled in such a way that indices increase by $1$ from one column to the next.

Now suppose we consider an arbitrary $2p \times p$ matrix $M$.  Let $\mathfrak{a}=(a_1,\dots,a_p)$ be a length $p$ sequence of row indices, and write $[\mathfrak{a}]_M$ (or simply $[\mathfrak{a}]$, if the underlying matrix $M$ is clear), for the  $p \times p$ minor obtained by selecting (in the corresponding order) the rows $a_1, \ldots, a_p$ 
of the matrix $M$. For instance, with $M=(x_{ij})_{i\leq 4,j\leq 2}$
we have
  $$[34] =\det\begin{bmatrix}
                      x_{{31}}&x_{{32}}\\
                      x_{{41}}&x_{{42}}
                    \end{bmatrix},
  \qquad {\rm and} \qquad 
      [32] =\det\begin{bmatrix}
                       x_{{31}}&x_{{32}}\\
                       x_{{21}}&x_{{22}}
                     \end{bmatrix}.$$
Then, for  any given $k$-subsequence\footnote{In the present context, such subsequences always inherit the order of the larger sequence they come from.} $\mathfrak{c}$ of $\mathfrak{a}$, $1\leq k\leq p$, the following {\em Pl\"ucker relation} holds on $p\times p$ minors of $M$:
\begin{equation}\label{plucker}
   [\mathfrak{a}]\,[\mathfrak{b}]=\sum_{\rouge{\mathfrak{d}}\subseteq_k \mathfrak{b}} \left[\mathfrak{a}|_{\bleu{\mathfrak{c}}\leftarrow \rouge{\mathfrak{d}}}\right]\,
                    \left[\mathfrak{b}|_{\rouge{\mathfrak{d}}\leftarrow  \bleu{\mathfrak{c}}}\right],
\end{equation}
with $\mathfrak{a}=(1,\ldots,p)$, $\mathfrak{b}=(p+1,\ldots, 2p)$, and the summation indices running over all $k$-subsequences $\mathfrak{d}$ of $\mathfrak{b}$. Here, we have used $\mathfrak{a}|_{\mathfrak{c}\leftarrow  \mathfrak{d}}$ to denote the sequence obtained by replacing in $\mathfrak{a}$ each index $c_i$, in $\mathfrak{c}$, by the index $d_i$, in $\mathfrak{d}$, in the position that $c_i$ appears in $\mathfrak{a}$.
For instance, with $n=4$ and $\mathfrak{c}=(\bf 1,3)$, we have the relation
\begin{eqnarray*}\bf
    [\bleu{1}2\bleu{3}4]\,[5678]&=&\bf
    [\rouge{5}2\rouge{6}4]\,[\bleu{13}78]+[\rouge{5}2\rouge{7}4]\,[\bleu{1}6\bleu{3}8]
    +[\rouge{5}2\rouge{8}4]\,[\bleu{1}67\bleu{3}]\\ 
          &&\bf\qquad + [\rouge{6}2\rouge{7}4]\,[5\bleu{13}8]
          +[\rouge{6}2\rouge{8}4]\,[5\bleu{1}7\bleu{3}]+[\rouge{7}2\rouge{8}4]\,[56\bleu{13}] .
\end{eqnarray*}
Before proving Proposition \ref{specialcase}, we wish to state and
prove a purely combinatorial lemma.
As usual, we define the {\em inversions} of a sequence
$\mathfrak{a}=(a_1, a_2, \ldots , a_p)$ of integers to be the pairs $(a_j, a_i)$ such that
$i<j$ and $a_i > a_j $. The number of inversions of $\mathfrak{a}$ is denoted ${\rm inv}(\mathfrak{a})$. 

\begin{lemma}\label{inversions}
Suppose we have two sequences $\mathfrak{a}=(a_1, \ldots, a_p)$ and 
$\mathfrak{b}= (b_1, \ldots, b_p)$ of integers satisfying
\begin{equation}\label{sequence}
a_1 \geq b_1 > a_2 \geq b_2 > \cdots > a_p \geq b_p.
\end{equation}
Pick any two length $k$ subsequences $\mathfrak{c}$ and $\mathfrak{d}$, respectively of $\mathfrak{a}$ and $\mathfrak{b}$. Then
$$
         {\rm inv}(\mathfrak{a}|_{\mathfrak{c}\leftarrow  \mathfrak{d}})=
         {\rm inv}(\mathfrak{b}|_{\mathfrak{d}\leftarrow  \mathfrak{c}}),
$$
if all entries of  $\mathfrak{a}|_{\mathfrak{c}\leftarrow  \mathfrak{d}}$ {\rm (}and of $\mathfrak{b}|_{\mathfrak{d}\leftarrow  \mathfrak{c}}${\rm )} are distinct.
\end{lemma}

\begin{proof}[\bf Proof.]
Assume that all the entries of 
   $$\mathfrak{a}'=(a'_1, \ldots, a'_p):=\mathfrak{a}|_{\mathfrak{c}\leftarrow  \mathfrak{d}}$$
    are distinct, as well as all the entries
of 
   $$\mathfrak{b}'= (b'_1, \ldots, b'_p):=\mathfrak{b}|_{\mathfrak{d}\leftarrow  \mathfrak{c}}.$$ 
Our first step is to reduce the problem to the case when all the inequalities
in \eqref{sequence} are strict.  Indeed, suppose $a_m=b_m$ for some
$m$.  Observe that this implies that $a_m$ and $b_m$ cannot be both
in $\mathfrak{a}'$ (or $\mathfrak{b}'$), since the elements of $\mathfrak{a}'$ are assumed
distinct. Then suppose we modify $\mathfrak{a}$ and $\mathfrak{b}$, and hence $\mathfrak{a}'$ and $\mathfrak{b}'$, 
by adding 1 to all
integers in \eqref{sequence} that are strictly greater 
than $a_m$, and then by adding 
1 to $a_m$ itself.  Thus, the only relative ordering in \eqref{sequence}
that is affected is that we now have $a_m > b_m$.  Since the original $a_m$ and $b_m$ could not both appear in the original $\mathfrak{a}'$, the number of inversions of the new $\mathfrak{a}'$
is unaffected by the modification (similarly for $\mathfrak{b}'$).  Therefore, we can assume
from now on that all the inequalities of \eqref{sequence} are strict.   Moreover, the condition on entries of $\mathfrak{a}'$ (or $\mathfrak{b}'$) being distinct is automatically satisfied in this case.

We prove the result by induction on $k$, with the result being trivially
true for $k=0$.  Let $a_\ell=c_k $ be the smallest integer in $\mathfrak{c}$,  and
$b_m=d_k $ be the smallest in $\mathfrak{d}$. In other words, these are the rightmost elements that are being switched.  Undoing this rightmost exchange, we get
 $${\bf a}:=\mathfrak{a}|_{\mathfrak{c}'\leftarrow  \mathfrak{d}'},\qquad {\rm  and} \qquad {\bf b}:=\mathfrak{b}|_{\mathfrak{d}'\leftarrow  \mathfrak{c}'},$$
 with $\mathfrak{c}'$ obtained by removing $a_\ell$ from $\mathfrak{c}$ (similarly for $\mathfrak{d}'$). Thus $a_\ell$ and $b_m$ are left in their original spots, and the pair $({\bf a},{\bf b})$ corresponds to an instance of the lemma with a smaller $k$. We need only consider he case when $\ell < m$, since  $\ell > m$ is similar. Clearly, the only inversions of $\mathfrak{a}'$ that are affected by the transition to ${\bf a}$ are those that involve $b_m$, who is sitting in position $\ell$. Among these, we can clearly restrict our considerations to entries of $\mathfrak{a}'$ with indices between $1$ and $m$; and we get exactly one such  inversion for each $x$ lying to the left of  $b_m$. Hence there are $\ell-1$ of these. By comparison, among the inversions in ${\bf a}$ arising from indices between $1$ and $m$, and involving $a_\ell$, we have  the following. Each entry between $1$ and $m$ ($\not= \ell$) gives rise to an inversion, except for the entries of $\mathfrak{d}$ with indices between $\ell$ and $m-1$. Say there are $k$ of these, then we get $m-1-k$ inversions involving $a_\ell$, among entries of ${\bf a}$ with index between $1$ and $m$. The difference in the number of inversion of ${\bf a}$ and $\mathfrak{a}'$ is thus $m-\ell-k$.
All this is illustrated below.
$$\begin{picture}(0,1)(12,-5)\linethickness{.8mm}
    \put(-2,-.2){$\mathfrak{a}'$}
    \put(11.5,.5){$\ell$}    \put(23.5,.5){$m$}
    \multiput(-.3,0)(3,0){9}{\jaune{\line(1,0){3}}}
    \multiput(-.3,0)(3,0){9}{\circle*{.3}}
    \put(2.65,0){\bleu{\circle*{.5}}}
    \put(8.65,0){\bleu{\circle*{.5}}}
    \put(11.65,0){\bleu{\circle*{.5}}}
    \put(3,-.4){\vector(3,-1){9}}
    \put(12,-3.4){\vector(-3,1){9}}
    \put(9,-.4){\vector(2,-1){6}}
    \put(15,-3.4){\vector(-2,1){6}}
    \put(12,-.4){\vector(4,-1){12}}
    \put(24,-3.4){\vector(-4,1){12}}
    \end{picture} 
    \begin{picture}(0,6)(11,-1)\linethickness{.8mm}
    \put(-2,-.2){$\mathfrak{b}'$}
    \put(11.5,-1.3){$\ell$}    \put(23.5,-1.2){$m$}
    \multiput(-.3,0)(3,0){9}{\jaune{\line(1,0){3}}}
    \multiput(-0.3,0)(3,0){9}{\circle*{.3}}
    \put(11.65,0){\rouge{\circle*{.5}}}
    \put(14.65,0){\rouge{\circle*{.5}}}
    \put(23.65,0){\rouge{\circle*{.5}}}
    \end{picture} 
    $$
A similar counting argument shows that the number of inversions of ${\bf b}$ differs from that of $\mathfrak{b}'$ by exactly the same quantity.  By the induction hypothesis, ${\bf a}$ and ${\bf b}$ have the same number of inversion, hence so have $\mathfrak{a}'$ and $\mathfrak{b}'$.\end{proof}

\begin{proof}[\bf Proof of Proposition \ref{specialcase}] Let $\gamma:=\mu\cup \nu$.
Without loss of generality, suppose that $\mu_1 = \ga_1$. We intend to consider  a Pl\"ucker relation of the form (\ref{plucker}) on $p\times p$ minors of  the following matrix:
 $$  H_\ga = \begin{bmatrix}
h_{\ga_1} & \ast & \cdots &\ast \\
\ast & h_{\ga_3} & \cdots & \ast \\
\vdots & \vdots & \ddots & \vdots \\
\ast & \ast & \cdots & h_{\ga_{2p-1}} \\ 
h_{\ga_2} & \ast & \cdots &\ast \\
\ast & h_{\ga_4} & \cdots & \ast \\
\vdots & \vdots & \ddots & \vdots \\
\ast & \ast& \cdots & h_{\ga_{2p}} 
\end{bmatrix}$$
where rows are completed just as in (\ref{JT_matrix}).
For this, we set $\mathfrak{c}$ to be the increasing sequence of indices in the set
  $$\left\{ 1 \leq i \leq p\ |\ \mu_i \neq \lt_i \right\},$$
say of cardinality $k$.
By the Jacobi-Trudi identity and definition (\ref{tilde_par}), the left-hand side of \eqref{plucker}
is simply
$s_\lt s_\mt$.  Consider the term of the right-hand side that occurs
when we choose the subsequence $(c_1+p,\ldots,c_k+p)$.  We see that this term is exactly $s_\mu s_\nu$.  Furthermore, all the other terms on the right-hand
side are clearly of the form $\pm s_\alpha s_\beta$, resulting in a formula taking the form
   \begin{equation}\label{expression}
     s_\lt s_\mt - s_\mu s_\nu=\sum_{\alpha,\beta} \pm s_\alpha s_\beta .
   \end{equation}
Each term in the right-hand side comes with a sign that depends on the order of
the rows in the minors considered. 
Since the product of two Schur functions is Schur-positive, it remains to
show that all these other terms appear with a plus sign, rather than
a minus sign.  

Just as they appear in (\ref{plucker}), we let $\mathfrak{a}:=(1,\ldots,p)$, $\mathfrak{b}:=(p+1,\ldots, 2p)$. For all $i$, $ 1 \leq i \leq 2p$, we define $\de(i)$ by
saying that the first term in row $i$ of the matrix $H_\ga$  is $h_{\de(i)}$.
More precisely, we set
    $$\de(i):=\begin{cases}
      \ga_{2i -1}-i+1,& \text{if\ }1\leq i\leq p \\
     \ga_{2(i-p)}-(i-p)+1, & \text{if\ } i>p .
\end{cases}$$
Also note that, since $\ga_1 \geq \ga_2 \geq \cdots \geq \ga_{2p}$, we
have that
  $$\de(1) \geq \de(p+1) > \de(2) \geq \de(p+2) > \cdots > \de(p) \geq \de(2p).$$
Now, for any length $k$ subsequence $\mathfrak{d}$ of $\mathfrak{b}$, we get 
    $$\left[\mathfrak{a}|_{{\mathfrak{c}}\leftarrow {\mathfrak{d}}}\right]= \pm s_\alpha,$$
for some partition $\alpha$, if all the rows in this minor are distinct.  
The zero minors have no impact on our discussion.
While we don't care specifically about $\alpha$, we need to know the sign preceding
$s_\alpha$. For any sequence $\mathfrak{c}=(c_1,\ldots, c_k)$, we consider the sequence
    $$\delta(\mathfrak{c}):=(\delta(c_1),\ldots, \delta(c_k)).$$
Then 
   $$\left[\mathfrak{a}'\right]= (-1)^{{\rm     
            inv}(\delta(\mathfrak{a}'))} s_\alpha,$$
with $\mathfrak{a}'=\mathfrak{a}|_{{\mathfrak{c}}\leftarrow {\mathfrak{d}}}$.
Similarly,
 $$\left[\mathfrak{b}'\right]= (-1)^{{\rm     
            inv}(\delta(\mathfrak{b}'))} 
s_\beta ,$$
with $\mathfrak{b}'=\mathfrak{b}|_{{\mathfrak{d}}\leftarrow {\mathfrak{c}}}$.
However, by Lemma \ref{inversions}, we know that 
 ${{\rm  inv}(\delta(\mathfrak{a}'))}={{\rm     inv}(\delta(\mathfrak{b}'))}$,
thus we conclude that each $s_\alpha s_\beta$ appears with a plus sign in (\ref{expression}).
\end{proof}


\section{Skew-Shape Generalization}\label{skewshape}

Inspired by a conjecture in \cite{BBR04pr}, one might ask if it makes sense to generalize
the $\sim$-operation to skew shapes.  
To ease our presentation, let us write
    $$(\alpha,\beta)\subseteq (\mu,\nu)\qquad {\rm whenever}\qquad
      \al \subseteq \mu \mbox{\ \ and \ } \be \subseteq \nu .$$
The following lemma is readily checked, and is
left as an exercise for the reader.  

\begin{lemma}\label{containment}
If $(\alpha,\beta)\subseteq (\mu,\nu)$, then $(\alpha,\beta)^\sim\subseteq (\mu,\nu)^\sim$.
\end{lemma}

A proof of this lemma is made easier if one considers the following equivalent {\em column definition}  of the $\sim$-operation.
Suppose $ (\lt, \mt)=(\mu, \nu) ^\sim$. 
Let us add zero parts $\mu'_i = 0$ to the conjugate $\mu'=(\mu_1',\mu_2',\ldots)$ of $\mu$, whenever $i$ is larger then the length $\ell(\mu')=\mu_1$ of $\mu'$. This will enable our statements to be length independent. 
Recall that $\mu'_i$ equals the number of rows of $\mu$ of length at least $i$. It follows easily from our original definition of the $\sim$-operation that
\begin{equation}\label{columndef}
\mu'_i + \nu'_i = \lt'_i + \mt'_i, \qquad{\rm and}\qquad \lt'_i - \mt'_i \in \{0,1\}.
\end{equation}
Furthermore, the conditions in \eqref{columndef} are sufficient to characterize $(\lt, \mt)$ uniquely.  In fact, we have
\begin{equation}\label{columndef2}
\lt'_i = \left\lceil \frac{\mu'_i + \nu'_i}{2} \right\rceil, \qquad{\rm and}\qquad \mt'_i = \left\lfloor \frac{\mu'_i + \nu'_i}{2} \right\rfloor.
\end{equation}
Intuitively, the $\sim$-operation simply has the effect of ``balancing out'' the 
column lengths, with a slight preference for $\lt$.  We will make much use
of the column definition 
in this section, where our main subject is evidence in
favor of the conjecture below.

Given skew shapes $\mu/\al$ and $\nu/\be$, 
Lemma \ref{containment} implies 
that it makes sense to set 
\[
(\mu/\al , \nu/\be)^\sim = (\lt/\at, \mt/\bt), 
\]
where $(\mu, \nu)^\sim = (\lt, \mt)$ and $(\al, \be)^\sim = (\at, \bt)$.
Considering the case when both $\al$ and $\be$ are the empty partition, 
we see that this is indeed a generalization of the $\sim$-operation for ordinary shapes.
Computer experiments, together with results presented in the remainder of this section, 
suggest that we state the following conjecture.

\begin{conjecture}\label{skewversion}
For all $\mu/\al$ and $\nu/\be$, if 
$(\mu/\al, \nu/\be)^\sim = (\lt/\at, \mt/\bt)$, then
\[
s_{\lt/\at} s_{\mt/\bt} - s_{\mu/\al} s_{\nu/\be}
\]
is Schur-positive.
\end{conjecture}

We should note that the skew shapes $\lt/\at$ and $\mt/\bt$ depend on the actual partitions $\mu$, $\nu$, $\al$
and $\be$ involved; and not just on the ``apparent''
skew shapes $\mu/\al$ and $\nu/\be$.  An example will be best to clarify this point.
Let us consider the two skew shapes
     $$21/1=\begin{picture}(2,2)(0,.5)
\thicklines
\put(1,0){\line(1,0){1}}
\put(0,1){\line(1,0){2}}
\multiput(1,1)(1,0){2}{\line(0,-1){1}}
\put(0,2){\line(1,0){1}}
\multiput(0,2)(1,0){2}{\line(0,-1){1}}
\thinlines
\put(0,0){\dashbox{0.2}(1,1)}
\end{picture}\,,\qquad
       32/2=\begin{picture}(3,2)(0,.5)
\thicklines
\put(2,0){\line(1,0){1}}
\put(1,1){\line(1,0){2}}
\multiput(2,1)(1,0){2}{\line(0,-1){1}}
\put(1,2){\line(1,0){1}}
\multiput(1,2)(1,0){2}{\line(0,-1){1}}
\thinlines
\put(0,0){\dashbox{0.2}(1,1)}
\put(1,0){\dashbox{0.2}(1,1)}
\put(0,1){\dashbox{0.2}(1,1)}
\end{picture} $$
which, up to translation, have the same configuration of boxes.
In the context of skew Schur functions, we actually have
   $$s_{21/1}=s_{32/2},$$
and it is usual to identify the two skew shapes, although they have different descriptions.  
By contrast, in our context, we may get different 
results from the skew version of the $\sim$-operation. For instance, we have 
\setlength{\unitlength}{3mm}
\begin{eqnarray*}
\left( 
\begin{picture}(2,2)(0,1)
\thicklines
\put(1,0){\line(1,0){1}}
\put(0,1){\line(1,0){2}}
\multiput(1,1)(1,0){2}{\line(0,-1){1}}
\put(0,2){\line(1,0){1}}
\multiput(0,2)(1,0){2}{\line(0,-1){1}}
\thinlines
\put(0,0){\dashbox{0.2}(1,1)}
\end{picture} 
\,,\ 
\begin{picture}(2,2)(0,1)
\thicklines
\put(0,0){\line(1,0){2}}
\put(0,1){\line(1,0){2}}
\multiput(0,1)(1,0){3}{\line(0,-1){1}}
\put(0,2){\line(1,0){1}}
\multiput(0,2)(1,0){2}{\line(0,-1){1}}
\end{picture} 
\right) ^\sim 
& = &
\left( 
\begin{picture}(2,2)(0,1)
\thicklines
\put(1,0){\line(1,0){1}}
\put(0,1){\line(1,0){2}}
\multiput(1,1)(1,0){2}{\line(0,-1){1}}
\put(0,2){\line(1,0){1}}
\multiput(0,2)(1,0){2}{\line(0,-1){1}}
\thinlines
\put(0,0){\dashbox{0.2}(1,1)}
\end{picture} 
\,,\ 
\begin{picture}(2,2)(0,1)
\thicklines
\put(0,0){\line(1,0){2}}
\put(0,1){\line(1,0){2}}
\multiput(0,1)(1,0){3}{\line(0,-1){1}}
\put(0,2){\line(1,0){1}}
\multiput(0,2)(1,0){2}{\line(0,-1){1}}
\end{picture} 
\right) , \\
\left( 
\begin{picture}(3,2)(0,1)
\thicklines
\put(2,0){\line(1,0){1}}
\put(1,1){\line(1,0){2}}
\multiput(2,1)(1,0){2}{\line(0,-1){1}}
\put(1,2){\line(1,0){1}}
\multiput(1,2)(1,0){2}{\line(0,-1){1}}
\thinlines
\put(0,0){\dashbox{0.2}(1,1)}
\put(1,0){\dashbox{0.2}(1,1)}
\put(0,1){\dashbox{0.2}(1,1)}
\end{picture} 
\ ,\ 
\begin{picture}(2,2)(0,1)
\thicklines
\put(0,0){\line(1,0){2}}
\put(0,1){\line(1,0){2}}
\multiput(0,1)(1,0){3}{\line(0,-1){1}}
\put(0,2){\line(1,0){1}}
\multiput(0,2)(1,0){2}{\line(0,-1){1}}
\end{picture} 
\right) ^\sim 
& = &
\left( 
\begin{picture}(3,2)(0,1)
\thicklines
\put(2,0){\line(1,0){1}}
\put(0,1){\line(1,0){3}}
\multiput(2,1)(1,0){2}{\line(0,-1){1}}
\put(0,2){\line(1,0){2}}
\multiput(0,2)(1,0){3}{\line(0,-1){1}}
\thinlines
\put(0,0){\dashbox{0.2}(1,1)}
\put(1,0){\dashbox{0.2}(1,1)}
\end{picture} 
\ ,\ 
\begin{picture}(2,2)(0,1)
\thicklines
\put(1,0){\line(1,0){1}}
\put(0,1){\line(1,0){2}}
\multiput(1,1)(1,0){2}{\line(0,-1){1}}
\put(0,2){\line(1,0){1}}
\multiput(0,2)(1,0){2}{\line(0,-1){1}}
\thinlines
\put(0,0){\dashbox{0.2}(1,1)}
\end{picture} 
\right)\,.
\end{eqnarray*}
We will say that $\mu/\al$ is a \emph{minimal pair} description 
of a skew shape  if both
\begin{itemize}\itemsep=4pt
\item[1)] $\al_i < \mu_i$, for all $1 \leq i \leq \ell(\al)$, and  
\item[2)] $\al'_j < \mu'_j$, for all $1 \leq j \leq \al_1$.
\end{itemize}
Thus, $21/1$ is a minimal pair, while
$32/21$ is not.  Using \cite{BucSoftware, SteSoftware}, we have verified Conjecture~\ref{skewversion} for all minimal pairs $\mu/\al$ and
$\nu/\be$, with $|\mu/\al| + |\nu/\be| \leq 12$.  This amounts to a total
of almost 1 million pairs of skew shapes.  
We emphasize, however, that
in Conjecture \ref{skewversion}, we do not require that
$\mu/\al$ and $\nu/\be$ be minimal.  

Recall that a \emph{horizontal strip} is a skew shape whose diagram has
at most one cell in each column.
Similarly, a \emph{vertical strip} has at most one cell is each row.
A \emph{ribbon} is a skew shape whose diagram is edgewise connected
and contains no $2 \times 2$ block of cells.  
Removing the connectedness restriction,
let us say that a \emph{weak ribbon}
is a skew shape whose diagram contains no $2 \times 2$ block of cells.  
Finally, a skew shape of the form $\mu/\al$, where $\mu$ and $\al$ are both non-empty
hooks, will be called a \emph{\skewedhook}. 

Notice that the $\sim$-operation does not preserve pairs of ribbons.
For example, 
\[
\left( 
\begin{picture}(2,2)(0,1)
\thicklines
\put(1,0){\line(1,0){1}}
\put(0,1){\line(1,0){2}}
\multiput(1,2)(1,0){2}{\line(0,-1){2}}
\put(0,2){\line(1,0){2}}
\put(0,2){\line(0,-1){1}}
\thinlines
\put(0,0){\dashbox{0.2}(1,1)}
\end{picture} 
\ ,\ 
\begin{picture}(1,1)(0,0.5)
\thicklines
\multiput(0,0)(1,0){2}{\line(0,1){1}}
\multiput(0,0)(0,1){2}{\line(1,0){1}}
\end{picture} 
\right) ^\sim 
= 
\left( 
\begin{picture}(2,2)(0,1)
\thicklines
\put(1,0){\line(1,0){1}}
\put(0,1){\line(1,0){2}}
\multiput(1,1)(1,0){2}{\line(0,-1){1}}
\put(0,2){\line(1,0){1}}
\multiput(0,2)(1,0){2}{\line(0,-1){1}}
\thinlines
\put(0,0){\dashbox{0.2}(1,1)}
\end{picture} 
\ ,\ 
\begin{picture}(2,1)(0,0.5)
\thicklines
\multiput(0,0)(0,1){2}{\line(1,0){2}}
\multiput(0,0)(1,0){3}{\line(0,1){1}}
\end{picture} 
\right).
\]
On the other hand, we have the following result.

\begin{proposition}\label{preserves}
The $\sim$-operation preserves the families of:
\begin{itemize}\itemsep=4pt
\item[(i)] pairs of horizontal strips,
\item[(ii)] pairs of vertical strips,
\item[(iii)] pairs of weak ribbons,
\item[(iv)] pairs of \skewedhook s.  
\end{itemize}
\end{proposition}

{\bf Proof.}  Throughout, let $(\mu/\al, \nu/\be)$ denote a pair 
of skew shapes of the designated form, and set $(\lt/\at, \mt/\bt)=(\mu/\al, \nu/\be)^\sim$.
\begin{itemize}\itemsep=4pt
\item[(i)]  We must show that $\lt'_i - \at'_i \leq 1$ and $\mt'_i-\bt'_i \leq 1$, for all $i$.  From
\eqref{columndef2}, we see that 
\[
\lt'_i = \left\lceil \frac{\mu'_i + \nu'_i}{2} \right\rceil \mbox{\ \ and\ \ }
\at'_i = \left\lceil \frac{\al'_i + \be'_i}{2} \right\rceil.
\]
But since $\mu/\al$ and $\nu/\be$ are horizontal strips, we know that 
$\mu'_i + \nu'_i - (\al'_i + \be'_i) \leq 2$ and so $\lt'_i - \at'_i \leq 1$.  Similarly,
$\mt'_i-\bt'_i \leq 1$.

\item[(ii)]  We must show that $\at'_i \geq \lt'_{i+1}$ for all $i$, and similarly for $\mt/\bt$.
Since $\mu/\al$ and $\nu/\be$ are vertical strips, we know that
$\al'_i \geq \mu'_{i+1}$ and $\be'_i \geq \nu'_{i+1}$.  The result now 
follows from \eqref{columndef2}.

\item[(iii)]  We must show that $\at'_i \geq \lt'_{i+1} - 1$, and similarly for $\mt/\bt$.
The argument is similar to that for (ii).

\item[(iv)]  This follows easily from the fact that the $\sim$-operation preserves pairs of hooks.
\qed
\end{itemize}

For any skew shape $\mu/\al$, let $\row{\mu/\al}$ (respectively $\col{\mu/\al}$) denote the partition whose multiset of parts equals the multiset of row (respectively column) lengths of $\mu/\al$.  In other words, if $\mu/\alpha$ is a minimal pair, the partition $\row{\mu/\alpha}$ is obtained by left justifying $\mu/\alpha$ on the $y$ axis and then reordering parts in decreasing order. This is illustrated in Figure \ref{row_just}.
\begin{figure}
    $${\rm row}\left(
\begin{picture}(4,3)(0,1.5)
\put(2,0){\line(1,0){2}}
\put(1,1){\line(1,0){3}}
\multiput(2,1)(1,0){3}{\line(0,-1){1}}
\put(0,2){\line(1,0){4}}
\multiput(1,2)(1,0){4}{\line(0,-1){1}}
\put(0,3){\line(1,0){2}}
\multiput(0,3)(1,0){3}{\line(0,-1){1}}
\put(0,4){\line(1,0){1}}
\multiput(0,4)(1,0){2}{\line(0,-1){1}}
\end{picture}\right)=
\begin{picture}(3,4)(0,1.5)
\put(0,0){\line(1,0){3}}
\put(0,1){\line(1,0){3}}
\multiput(0,1)(1,0){4}{\line(0,-1){1}}
\put(0,2){\line(1,0){2}}
\multiput(0,2)(1,0){3}{\line(0,-1){1}}
\put(0,3){\line(1,0){2}}
\multiput(0,3)(1,0){3}{\line(0,-1){1}}
\put(0,4){\line(1,0){1}}
\multiput(0,4)(1,0){2}{\line(0,-1){1}}
\end{picture}\,,\qquad 
{\rm col}\left(
\begin{picture}(4,3)(0,1.5)
\put(2,0){\line(1,0){2}}
\put(1,1){\line(1,0){3}}
\multiput(2,1)(1,0){3}{\line(0,-1){1}}
\put(0,2){\line(1,0){4}}
\multiput(1,2)(1,0){4}{\line(0,-1){1}}
\put(0,3){\line(1,0){2}}
\multiput(0,3)(1,0){3}{\line(0,-1){1}}
\put(0,4){\line(1,0){1}}
\multiput(0,4)(1,0){2}{\line(0,-1){1}}
\end{picture}\right)=
\begin{picture}(4,2)(0,.5)
\put(0,0){\line(1,0){4}}
\put(0,1){\line(1,0){4}}
\multiput(0,1)(1,0){5}{\line(0,-1){1}}
\put(0,2){\line(1,0){4}}
\multiput(0,2)(1,0){5}{\line(0,-1){1}}
\end{picture}\,.$$
 \caption{Row justification versus column justification.}\label{row_just}
 \end{figure}
The following result is central to the proofs of this section.  

\begin{lemma}\label{skewdominance}
Suppose $(\mu/\al, \nu/\be)^\sim = (\lt/\at, \mt/\bt)$.  Then 
\begin{itemize}\itemsep=4pt
\item[(i)] $\row{\mu/\al} \cup \row{\nu/\be} \domgeq \row{\lt/\at} \cup \row{\mt/\bt}$,
\item[(ii)] $\col{\mu/\al} \cup \col{\nu/\be} \domgeq \col{\lt/\at} \cup \col{\mt/\bt}$.
\end{itemize}
\end{lemma}

We postpone the proof until the end of this section.  This lemma has a number of 
important implications, as we now begin to explain.

As for ordinary shapes, let us define the \emph{support}, $\supp{\mu/\al}{\nu/\be}$, of 
$\mu/\al$ and $\nu/\be$ to be the set of all partitions $\theta$ such that $s_\theta$ 
appears, with non-zero coefficient, in the product $s_{\mu/\al}s_{\nu/\be}$.  
It is clear that a necessary condition for Conjecture~\ref{skewversion} to be true
is that $\supp{\mu/\al}{\nu/\be} \subseteq \supp{\lt/\at}{\mt/\bt}$.  The following
corollary of Lemma~\ref{skewdominance} says that the extreme elements
of the supports are consistent with Conjecture~\ref{skewversion}. 

\begin{corollary}
For a pair of skew shapes $(\mu/\al, \nu/\be)$, the set $\supp{\mu/\al}{\nu/\be}$
has a unique maximum element and a unique minimum element in dominance order,
which we denote respectively by $\max{\supp{\mu/\al}{\nu/\be}}$
and $\min{\supp{\mu/\al}{\nu/\be}}$.  
If 
$(\lt/\at, \mt/\bt)=(\mu/\al, \nu/\be)^\sim$ then we also have:
\begin{itemize}\itemsep=4pt
\item[(i)] $\max{\supp{\mu/\al}{\nu/\be}} \domleq \max{\supp{\lt/\at}{\mt/\bt}}$,
\item[(ii)] $\min{\supp{\mu/\al}{\nu/\be}} \domgeq \min{\supp{\lt/\at}{\mt/\bt}}$.
\end{itemize}
\end{corollary}

\begin{proof}[\bf Proof.]
Consider the filling of $(\mu/\al) * (\nu/\be)$ that results from filling the $i$th lowest
cell of each column with the number $i$.  It is not difficult to see that this gives
a semistandard tableau which is an LR-filling.  It follows that 
$\pi := (\col{\mu/\al} \cup \col{\nu/\be})' \in \supp{\mu/\al}{\nu/\be}$.  Now consider
any $\theta \in \supp{\mu/\al}{\nu/\be}$ which results from an LR-filling $t$ of
$(\mu/\al) * (\nu/\be)$.  We have $\theta_1 \leq \pi_1$, since $t$ can have at most one 1
in each column.  In general, $\theta_1 + \cdots +\theta_i \leq \pi_1 + \cdots +\pi_i$,
since $t$ can have at most $i$ entries less than or equal to $i$ in each column.  
We conclude that $\max{\supp{\mu/\al}{\nu/\be}}$ exists and equals 
$(\col{\mu/\al} \cup \col{\nu/\be})'$.  Applying Lemma~\ref{skewdominance}(ii), we
get that $\max{\supp{\mu/\al}{\nu/\be}} \domleq \max{\supp{\lt/\at}{\mt/\bt}}$, proving (ii).

To prove (i), we exploit the proof of (ii).  Suppose $\theta \in \supp{\mu/\al}{\nu/\be}$.  
Applying the involution $\omega$, we see that $\theta \in \supp{\mu/\al}{\nu/\be}$
if and only if $\theta' \in \supp{\mu'/\al'}{\nu'/\be'}$.  It follows that
\[
\theta' \domleq (\col{\mu'/\al'} \cup \col{\nu'/\be'})'
\]
and so 
\begin{eqnarray*}
\theta & \domgeq & \col{\mu'/\al'} \cup \col{\nu'/\be'} \\
& = & \row{\mu/\al} \cup \row{\nu/\be}.
\end{eqnarray*}
Therefore, $\min{\supp{\mu/\al}{\nu/\be}}$ exists and equals $\row{\mu/\al} \cup \row{\nu/\be}$.
Applying Lemma \ref{skewdominance}(i), we conclude (ii).

\end{proof}

\begin{theorem}\label{skewcases}
Conjecture~\ref{skewversion} holds when $(\mu/\al, \nu/\be)$ is 
\begin{itemize}\itemsep=4pt
\item[(i)] a pair of horizontal strips,
\item[(ii)] a pair of vertical strips,
\item[(iii)] a pair of \skewedhook s.
\end{itemize}
\end{theorem}

\begin{proof}[\bf Proof.] \begin{enumerate}\itemsep=4pt
\item[(i)] We first observe that, when $\mu/\alpha$ is an horizontal strip, we have 
    $$s_{\mu/\alpha}=h_{\row{\mu/\alpha}}.$$
Thus, when $(\mu/\al, \nu/\be)$ is a pair of horizontal strips, $s_{\mu/\al}s_{\nu/\be}$
depends only on 
     $$\row{\mu/\al} \cup \row{\nu/\be}=\row{(\mu/\alpha) *(\nu/\beta)}.$$
By Proposition~\ref{preserves}(i), $ (\lt/\at, \mt/\bt)=(\mu/\al, \nu/\be)^\sim $ is also 
a pair of horizontal strips.  Therefore, 
\[
s_{\lt/\at} s_{\mt/\bt} - s_{\mu/\al} s_{\nu/\be} = 
h_{\row{\lt/\at} \cup \row{\mt/\bt}} - h_{\row{\mu/\al} \cup \row{\nu/\be}}.
\]
Now,  $h_\theta - h_\pi$ is Schur-positive
if and only if $\theta \domleq \pi$ (see, for example, \cite[p. 119]{Mac95}).
Therefore, (i) follows from Lemma~\ref{skewdominance}(i).

\item[(ii)] The proof is similar to that of (i), expect that we now use the fact that for a vertical strip $\mu/\alpha$, $s_{\mu/\alpha} = e_{{\rm col}(\mu/\alpha)}$.
Also, applying the involution $\omega$, we see that the Schur-positivity of 
$e_\theta - e_\pi$ is equivalent to that of $h_\theta - h_\pi$.  Given these
two facts, we can use the argument of (i) to deduce (ii) from 
Lemma~\ref{skewdominance}(ii).

\item[(iii)]  
Suppose $\mu=(M,1^m)$, $\nu=(N,1^n)$, $\al=(A, 1^a)$ and $\be=(B, 1^b)$.  Without
loss of generality, suppose $M \geq N$.  Let us set 
$$S:=\max\{A,B\},\qquad
T:=\min\{A,B\},$$
and
 $$\ell:=\left\lceil \frac{m+n}{2} \right\rceil,\quad 
      r:=\left\lfloor \frac{m+n}{2} \right\rfloor,\quad
        s:=\left\lceil \frac{a+b}{2} \right\rceil,\quad
        t:=\left\lfloor \frac{a+b}{2} \right\rfloor .$$
We easily see, from the definition, that we have $(\lt/\at, \mt/\bt)=(\mu/\al,\nu/\be)^\sim $, with
\begin{equation*}
\lt = (M, 1^\ell), \quad
\mt = (N, 1^r), \quad
\at = (S, 1^s), \quad
\bt= (T, 1^t).
\end{equation*}
Since $\al$ and $\be$ are non-empty by definition of \skewedhook s, we get
\begin{eqnarray*}
s_{\lt/\at} s_{\mt/\bt} - s_{\mu/\al} s_{\nu/\be}
 & = & 
h_{M-S}h_{N-T}e_{\ell-s} 
e_{r-t} - h_{M-A}h_{N-B}e_{m-a}e_{n-b} \\
& = & 
h_{M-S}h_{N-T} (e_{\ell-s} 
e_{r-t}  - e_{m-a}e_{n-b}) \\
& &\quad + e_{m-a}e_{n-b}(h_{M-S}h_{N-T} - h_{M-A}h_{N-B}).
\end{eqnarray*}
Now both terms in parentheses in this latter expression are Schur-positive.  Indeed, 
consider the pair of vertical strips $((1^m)/(1^a), (1^n)/(1^b))$.  From (ii), we deduce
that $e_{\ell-s} e_{r-t} - e_{m-a}e_{n-b}$ is Schur-positive.
Similarly, considering $((M)/(A), (N)/(B))$, we deduce from (i) 
that $h_{M-S}h_{N-T} - h_{M-A}h_{N-B}$ is Schur-positive.  Thus
$s_{\lt/\at} s_{\mt/\bt} - s_{\mu/\al} s_{\nu/\be}$ is Schur-positive, as required.
\end{enumerate}\end{proof}

Before concluding this section with the proof of Lemma~\ref{skewdominance}, 
we need two facts about the dominance order.  We first introduce some notation.
We extend the addition notation used for partitions to sequences of integers of the same length. We also observe that the definition of the dominance order,
exactly as it is as stated in (\ref{def_dom}),
can obviously be extended to weakly decreasing sequences of integers.  
For a sequence of integers $C$, let $\partition{C}$ denote the sequence obtained by sorting the entries of $C$ into weakly decreasing order. In particular, if $C$ consists of non-negative integers,
we get a partition.

\begin{lemma}\label{dominancepreserved}
\begin{itemize}\itemsep=4pt
\item[(i)] If $\mu$ and $\nu$ are partitions with $\mu \domgeq \nu$, then $\addn{\mu}{d} \domgeq 
\addn{\nu}{d}$.
\item[(ii)] Let $\gamma = (\gamma_1, \ldots, \gamma_\ell)$ be a weakly decreasing sequence
of integers, and $\delta = (\delta_1, \ldots, \delta_\ell)$ be a weakly increasing sequence of 
integers. Then, for any permutation  $\varepsilon$ of the sequence $\delta$, we have
 $\partition{(\gamma + \varepsilon)} \domgeq \partition{(\gamma + \delta)}$.
\end{itemize}
\end{lemma}

\begin{proof}[\bf Proof.]
(i) Suppose $(\addn{\mu}{d})_r = d$ while $(\addn{\nu}{d})_s = d$.  The case $r=s$ is simple to show.
Suppose $r > s$.  We see easily that
\[
\sum_{i=1}^{j} (\addn{\mu}{d})_i \geq \sum_{i=1}^{j} (\addn{\nu}{d})_i
\]
when $j < s$ or when $j \geq r$.  Therefore, assume $s \leq j < r$.  We have
\begin{eqnarray*}
\sum_{i=1}^{j} (\addn{\mu}{d})_i & \geq & \sum_{i=1}^{s-1} (\addn{\mu}{d})_i + (j-(s-1))d \\
& \geq & \sum_{i=1}^{s-1} (\addn{\nu}{d})_i +(j-(s-1))d \\
& \geq & \sum_{i=1}^{j} (\addn{\nu}{d})_i
\end{eqnarray*}
as required.  
The case $r < s$ is somewhat similar, yet the main idea is different enough to warrant a 
demonstration.  We know that 
\[
\sum_{i=1}^{j} (\addn{\mu}{d})_i \geq \sum_{i=1}^{j} (\addn{\nu}{d})_i
\]
when $j < r$ or when $j \geq s$.  Therefore, assume $r \leq j < s$.  We have
\begin{eqnarray*}
\sum_{i=1}^{j} (\addn{\mu}{d})_i & \geq & \sum_{i=1}^{s} (\addn{\mu}{d})_i - (s-j)d \\
& \geq & \sum_{i=1}^{s} (\addn{\nu}{d})_i - (s-j)d \\
& \geq & \sum_{i=1}^{j} (\addn{\nu}{d})_i 
\end{eqnarray*}
as required.

(ii)  We proceed by induction on $\ell$, with the result being trivially true for $\ell = 1$.
Suppose $\varepsilon_1 = \delta_s$, and define $\zeta = (\delta_s, \delta_1, \ldots, \delta_{s-1}, \delta_{s+1}, \ldots, \delta_\ell)$.  Our approach is to show that
\[
\partition{(\gamma+\varepsilon)} \domgeq \partition{(\gamma+\zeta)} \domgeq
\partition{(\gamma+\delta)}.
\]
We first show that $\partition{(\gamma+\varepsilon)} \domgeq \partition{(\gamma+\zeta)}$.  
Suppose that $\gamma_1 + \varepsilon_1 = \gamma_1 + \zeta_1 = d$.  Let
$\bar{\gamma} := (\gamma_2, \ldots, \gamma_\ell)$, with 
$\bar{\varepsilon}$ and $\bar{\zeta}$ defined similarly.
Now $\bar{\zeta}$ is weakly decreasing and $\bar{\varepsilon}$
is a permutation of $\bar{\zeta}$, so by the induction hypothesis, 
$\partition{(\bar{\gamma} + \bar{\varepsilon})} \domgeq \partition{(\bar{\gamma} + \bar{\zeta})}$.
Applying part (i) of this lemma, we deduce that
$\partition{(\gamma+\varepsilon)} \domgeq \partition{(\gamma+\zeta)}$.

It remains to show that $\partition{(\gamma+\zeta)} \domgeq \partition{(\gamma+\delta)}$.
If there exists $i \leq \ell$ such that $\zeta_i = \delta_i$, then we can deduce the result
by applying the same argument as in the previous paragraph.  Assume, therefore, 
that $s=\ell$ and so $\zeta = (\delta_\ell, \delta_1, \ldots, \delta_{\ell-1})$.  

Supposing that the $j$ largest elements of the sequence $\gamma + \delta$ are in 
positions 
    $$i_1 < i_2 < \cdots < i_j,$$
we have
\begin{eqnarray*}
\sum_{i=1}^{j} (\partition{(\gamma + \zeta)})_i  & \geq &
(\gamma_1 + \delta_\ell) + (\gamma_{i_2} + \delta_{i_2 -1}) + (\gamma_{i_3} + \delta_{i_3 -1})  
+ \cdots + (\gamma_{i_j} + \delta_{i_j -1}) \\
& = & (\gamma_{i_1}+\delta_{i_1}) + (\gamma_{i_2}+\delta_{i_2}) + \cdots (\gamma_{i_j}+\delta_{i_j}) \\
& & + (\gamma_1 - \gamma_{i_1}) + (\delta_{i_2 -1}-\delta_{i_1}) 
+ (\delta_{i_3 -1} - \delta_{i_2}) + \cdots + (\delta_{i_j - 1} - \delta_{i_{j-1}}) \\
& & + (\delta_\ell - \delta_{i_j}) \\
& \geq & (\gamma_{i_1}+\delta_{i_1}) + (\gamma_{i_2}+\delta_{i_2}) + 
\cdots + (\gamma_{i_j}+\delta_{i_j}) \\
& = & \sum_{i=1}^{j} (\partition{(\gamma+\delta)})_i \ .
\end{eqnarray*}
\end{proof}

\begin{proof}[\bf Proof of Lemma~\ref{skewdominance}.]

(i) Let $\gamma = \mu \cup \nu = \lt \cup \mt$, which we suppose to be of length $\ell$.
Let $\delta$ be the weakly increasing
sequence of integers with elements consisting of the multiset
$\{-(\al \cup \be)_i \}_{i=1}^{\ell}$.  The reason for this somewhat contrived definition
of $\delta$ is we get that $\gamma+\delta = \row{\lt/\at} \cup \row{\mt/\bt}$.  
Furthermore, $\row{\mu/\al} \cup \row{\nu/\be} = \gamma+\varepsilon$, where 
$\varepsilon$ is some permutation of $\delta$.  We see that we are in exactly the 
situation necessary to apply Lemma~\ref{dominancepreserved}(ii).  It
follows that $\row{\mu/\al} \cup \row{\nu/\be} \domgeq \row{\lt/\at} \cup \row{\mt/\bt}$,
as required.

(ii) Consider the $i$th columns of $\mu/\al$ and $\nu/\be$.  They contribute the
multiset $A_i = \{\mu'_i-\al'_i, \nu'_i-\be'_i\}$ to $\col{\mu/\al} \cup \col{\nu/\be}$.  
By \eqref{columndef2}, the $i$th columns of $\lt/\at$ and $\mt/\bt$ contribute
the multiset $B_i$ to $\col{\lt/\at} \cup \col{\mt/\bt}$, where
\[
B_i = \left\{ \left\lceil \frac{\mu_i + \nu_i}{2} \right\rceil - 
\left\lceil \frac{\al_i + \be_i}{2} \right\rceil  , 
\left\lfloor \frac{\mu_i + \nu_i}{2} \right\rfloor -
\left\lfloor \frac{\al_i + \be_i}{2} \right\rfloor
\right\}.
\]
We claim that 
$\partition{A_i} \domgeq \partition{B_i}$.  Indeed, $|\partition{A_i}| = |\partition{B_i}|$.
Therefore, it suffices to observe that $(\partition{B_i})_1 - (\partition{B_i})_2 \in \{0, 1\}$.

By repeated applications of Lemma~\ref{dominancepreserved}(i), we see that for any partitions
$\lambda$, $\al$ and $\be$, if $\al \domgeq \be$, then $\al \cup \lambda \domgeq \be \cup \lambda$.  
Furthermore, for any $\at \domgeq \bt$,
\[
\al \cup \at \domgeq \al \cup \bt \domgeq \be \cup \bt .
\]
Repeatedly applying this to our partitions $\partition{A_i}$ and $\partition{B_i}$, we conclude that
\[
\col{\mu/\al} \cup \col{\nu/\be} = \bigcup_i \partition{A_i} \domgeq 
\bigcup_i \partition{B_i} = \col{\lt/\at} \cup \col{\mt/\bt} .
\]
\end{proof}

\begin{remark}\label{skewmtilde}
As mentioned in the introduction, all the results of this section hold when we extend from
pairs of partitions to $m$-tuples of partitions, for $m \geq 2$.  There are two ways
to see this.  The first is to observe
that all our proofs of results specific to the $m=2$ case also 
work in the extended case after minor modifications.  
Alternatively, we can note that the repainting argument referred to in the introduction
only modifies a pair of partitions at each stage.  Because of this, it can be used
to extend Lemma~\ref{containment}, Proposition~\ref{preserves} and Theorem~\ref{skewcases}
to the general $m$ case.  
\end{remark}


\section{Posets of Pairs of Partitions}\label{posets}

\begin{figure}
$$\scalebox{.6}{\includegraphics{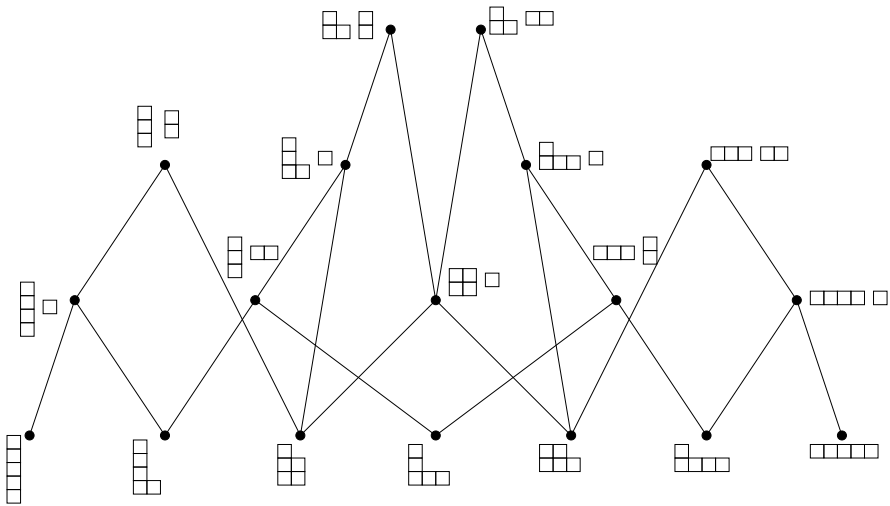}}$$
\caption{$\p{5}$.  Note that empty partitions are not displayed.}
\label{allpairs5}
\end{figure}

In studying the $\ast$-operation and the $\sim$-operation, the main difficulty lies
in  understanding
when an expression of form \eqref{diff_pairs} is Schur-positive.  To take a global 
approach to this question, fix $n$ and consider the set $\p{n}$ of all pairs of partitions $(\mu, \nu)$
such that $|\mu| + |\nu| = n$.  For the purposes of the current discussion, we identify 
$(\mu, \nu)$ and $(\nu, \mu)$.  We make $\p{n}$ into a poset by saying that $(\mu, \nu) \leq (\tau, \sigma)$
if  $s_\tau s_\sigma - s_\mu s_\nu$ is Schur-positive.  It is possible, but not necessarily easy, 
to see that $\leq$ is then an antisymmetric relation.  Figure \ref{allpairs5} shows $\p{5}$.  

While understanding this poset for general $n$ may seem like a difficult task, it is encouraging 
to look at the following similarly defined poset: say that $(\mu, \nu) \leq (\tau, \sigma)$
if $h_\tau h_\sigma - h_\mu h_\nu$ is Schur-positive, where $h$ denotes the complete
homogeneous symmetric functions.  To maintain the poset structure, we identify all those
pairs $(\mu, \nu)$ and $(\tau, \sigma)$ such that $\mu \cup \nu = \tau \cup \sigma$, since in that case we evidently have  $h_\tau h_\sigma - h_\mu h_\nu = 0$.  As we earlier observed, 
$h_\al - h_\be$ is Schur-positive if and only if $\al \domleq \be$.
Hence the resulting poset is exactly the (self-dual) dominance lattice for partitions of $n$. 

\begin{figure}
\center
$$\scalebox{.6}{\includegraphics{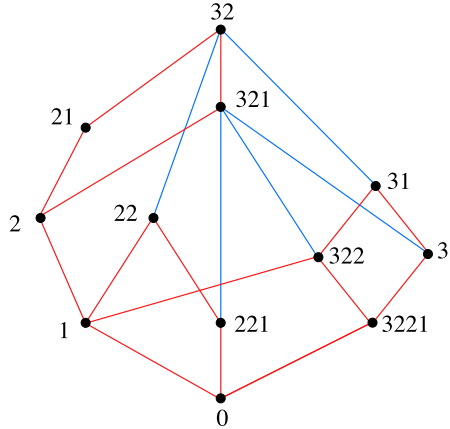}}$$
\caption{The poset $\pd{k3221}$, for $k\geq 5$.}
\label{pd53221}
\end{figure}

It is insightful to consider Conjecture~\ref{fflp27} in terms of the posets $\p{n}$.
Given a partition $\ga=(\ga_1, \ga_2, \cdots, \ga_{2p})$,
with $|\ga| = n$ and allowing $0$ parts, we wish to know which ``dealing'' of the parts of $\ga$ between
two partitions $\mu$ and $\nu$ will result in a large pair $(\mu, \nu)$ in the
poset $\p{n}$.  
Let $\pd{\ga}$ denote the subposet of $\p{n}$ consisting of
those pairs $(\mu, \nu)$ that arise as a dealing of the parts of $\ga$.  Then
Conjecture~\ref{fflp27} states that $\pd{\ga}$ has a unique maximal 
element for every $\ga$, namely $(\lt, \mt)$, where
\[
\lt = (\ga_1, \ga_3, \ldots , \ga_{2p-1}), \ \ \ 
\mt=(\ga_2, \ga_4, \ldots , \ga_{2p}).
\]
Figure \ref{pd53221} illustrates the poset  $\pd{53221}$.  We label each element $(\mu, \nu)$ by $\nu$ 
only, where without loss of generality, we always take $\nu$ to be the partition which is smaller
in lexicographic order ($3221$ rather then $5$), with no regard for the parts sum. 

\begin{figure}
\center
$$\scalebox{.6}{\includegraphics{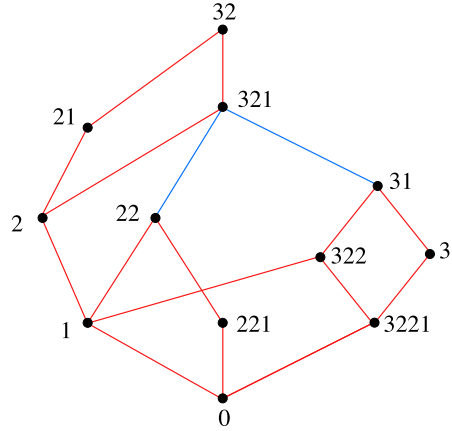}}$$
\caption{The poset $\pd{43221}$.}
\label{pd43221}
\end{figure}

Interestingly, our investigations suggest that, in a large number of cases, $\pd{\ga}$
is independent of the sizes of the parts of $\ga$ and only depends on the set of indices
$i$ such that $\ga_i$ is \emph{strictly} greater than $\ga_{i+1}$.  For example, there
is a canonical isomorphism from $\pd{53221}$ to $\pd{k3221}$ for any $k \geq 5$. Actually, using our labeling convention above, Figure \ref{pd53221} gives the poset associated to any partition of the form $k3221$ with $k\geq 5$.
This will be a consequence of Proposition \ref{pdinvariance} below.  On the other 
hand,  $\pd{53221}$ is not isomorphic to the poset $\pd{43221}$, which is illustrated in Figure \ref{pd43221}.  Rather, as is readily seen, it is a weak subposet: i.e.
the set of order relations in $\pd{53221}$ is a subset of the set of relations in $\pd{43221}$.

\begin{proposition}\label{pdinvariance}
Let $\ga = (\ga_1, \ga_2, \ldots , \ga_m)$ be a partition with $\ga_1 \geq \ga_2 + \ga_3$.
Then  $$\pd{\ga} \cong \pd{\gamma_1+1, \ga_2, \ldots, \ga_m}.$$
\end{proposition}

\begin{proof}[\bf Proof.]
Let $(\lambda, \mu)$ be an element of $\pd{\ga}$ where, without loss of generality, we take
$\lambda_1 \geq \mu_1$.  For any partition $\nu$, let $\addone{\nu}$ denote the partition 
obtained when the first part of $\nu$ is increased by 1.
Consider those partitions $\theta$ such that $c^\theta_{\addone{\lambda} \mu} > 0$.
We have that 
 $$\theta_1 \geq \lambda_1 + 1 = \ga_1 + 1 \geq \ga_2 + \ga_3 +1 \geq \mu_1 + \lambda_2 + 1 > \theta_2,$$
 where the last inequality can be seen by considering LR-fillings of 
 $\theta/\mu$ of content $\lambda$.
Therefore, $\theta = \addone{\nu}$ for some $\nu$, and we write $\nu=(\theta-1)$ when this is the case. We claim that, for any $\theta$ such that $c^\theta_{\addone{\lambda} \mu} > 0$ as above, there is a bijection:
 $$ \mbox{LR-fillings of $(\theta-1)/\mu$ of content $\lambda$}
 \longleftrightarrow
 \mbox{LR-fillings of $\theta/\mu$ of content $\addone{\lambda}$} . $$
We see that the existence of such a bijection, $\varphi$, would imply the desired canonical poset isomorphism.
Given any LR-filling $t$, of shape $(\theta-1)/\mu$ and content $\lambda$, let $\varphi(t)$ be the semistandard  tableau of $\theta/\mu$ obtained from $t$ simply
by adding a $1$ in the first row to the right of the last entry of $t$.  We easily see that this is indeed an LR-filling of $\theta/\mu$, and has content $\addone{\lambda}$ by construction.  We must check that  the inverse map makes sense and  has the required properties.  
If $u$ is an LR-filling of $\theta/\mu$ of content $\addone{\lambda}$, then the inverse of $\varphi$ corresponds to the deletion of the rightmost $1$ in the first row of $u$.  We see that this gives a semistandard tableau $\varphi^{-1}(u)$, with content $\lambda$ by construction. Hence, it only remains to check that $\varphi^{-1}(u)$ is actually an LR-filling.
The number of $1$'s in the first row of $\varphi^{-1}(u)$ equals $(\theta_1-1) - \mu_1$.  
If $\mu_1 = \ga_2$ then we have 
  $$(\theta_1-1) - \mu_1 \geq \lambda_1 - \ga_2 = \ga_1 - \ga_2 \geq \ga_3 \geq \lambda_2.$$
Otherwise, $\mu_1 \leq \ga_3$ and we get
  $$(\theta_1-1) - \mu_1 \geq \lambda_1 - \ga_3 = \ga_1 - \ga_3 \geq \ga_2 = \lambda_2.$$
In either case, the number of $1$'s in the first row of $\varphi^{-1}(u)$ is greater than or equal to the  total number
of $2$'s in $\varphi^{-1}(u)$, as required.  Therefore, $\varphi^{-1}(u)$ is an LR-filling, and this finishes the proof.
\end{proof}

\begin{figure}
$$\scalebox{.6}{\includegraphics{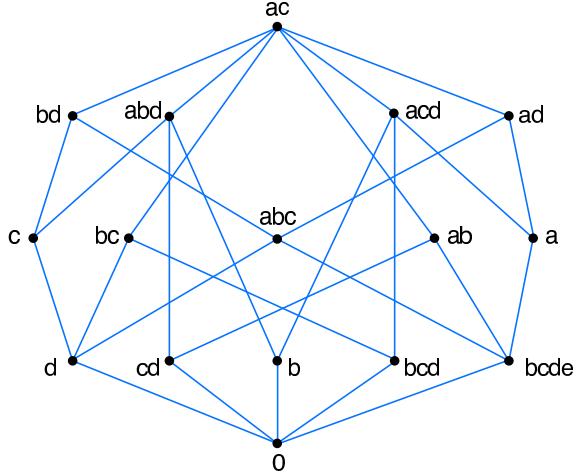}}$$
\caption{Poset $\pd{x,a,b,c,d}$, when $10\geq x>a>b>c>d>0$.}
\label{abcd}
\end{figure}

There seems to be a natural chamber complex decomposition, of the space of partitions of length $k$, induced by isomorphism classes of associated posets.  Proposition \ref{pdinvariance} is an example of results along these lines, since it states that        
           $$\gamma_1=\gamma_2+\gamma_3$$
might be
one of the defining hyperplanes of such a chamber complex. Various computational experiments also confirm this impression. For instance, for all partition of length $5$, into distinct parts (all $\leq 10$), we always get posets isomorphic to that of Figure \ref{abcd}.


\section{An Explicit Example with an Application}\label{explodedjt}

One interesting special case for which this general invariance, with respect to relative part size,  of the poset $\pd{\ga}$ holds, is the following. As we will see, even this special case implies nice identities.
\begin{proposition}
For all $k$ and $m$, $\pd{k^m}$ is isomorphic to a chain of length  $\ell= \lfloor \frac{m}{2} \rfloor$.  
\end{proposition}

\begin{proof}[\bf Proof.]
The elements of  $\pd{k^m}$ are of the form  $d_r:=((k^r),(k^{m-r}))$, with
$r \leq \ell$, and we need to show that $d_{r-1} < d_r$ for $r = 1, \ldots, \ell$.
To show the proposition, we need only construct, for each partition $\theta$ such that $c^\theta_{(k^{m-r+1}), (k^{r-1})} > 0$, an injection
    $$\psi:\mathcal{L}_{r-1} \hookrightarrow \mathcal{L}_{r},$$ 
 from the set $\mathcal{L}_{r-1}$, of LR-fillings of $\theta / (k^{m-r+1})$ with content $(k^{r-1})$, to  the set $\mathcal{L}_{r}$ of LR-fillings of $\theta/ (k^{m-r})$ with content $(k^r)$. By definition, column entries of any semistandard skew tableau of shape $\theta/ (k^{m-r})$ have to be strictly increasing. Thus, to have content $(k^{r-1})$, the maximal height of its columns has to be bounded by $r-1$. Moreover, any prefix of an LR-filling, is an LR-filling.  Since $m-r \geq r$, it follows that an LR-filling $t$ in $\mathcal{L}_{r-1}$ consists of two disjoint ``pieces'': an LR-filling 
of the shape $\alpha$ sitting to the right of the rectangle $(k^{m-r+1})$, and a filling 
of the shape $\beta$ sitting above this same rectangle.  In view of the preceding argument,
the number of parts of both $\alpha$ and $\beta$ are at most $r-1$, . In particular, this also implies that $\beta_1 \leq k$. The typical configuration is illustrated in Figure \ref{lrfilling}(a).

Now, we construct $\psi(t)$ as follows. First, convert the top (length $k$) row of 
$( k^{m-r+1})$ to a row $r$'s and reorder each column of the
resulting tableau so 
that it is increasing as we go up. The result is still
a semistandard tableau. This process is illustrated Figure \ref{lrfilling}(b).
\setlength{\unitlength}{4mm}\begin{figure}
\center
%
$$\hspace{-8mm} \scalebox{.8}{\includegraphics{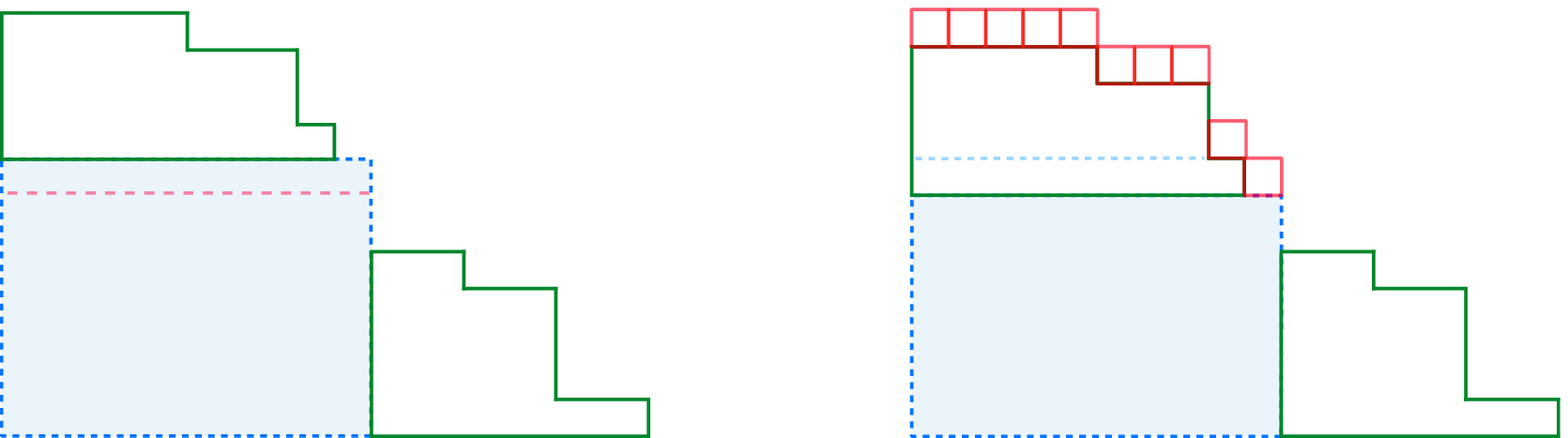}}$$
\caption{}
\label{lrfilling}
\begin{picture}(0,0)(0,0)
\put(-14,11){$\beta$}
\put(5.7,10.3){$\beta$}
\put(-7,5.5){$\alpha$}
\put(12.7,5.5){$\alpha$}
\put(-15,7){$(k^{m-r+1})$}
\put(5,7){$(k^{m-r})$}
\multiput(3,12.6)(.8,0){5}{$r$}
\multiput(7,11.8)(.8,0){3}{$r$}
\put(9.4,10.2){$r$}\put(10.2,9.4){$r$}
\put(-12,2.5){(a)}\put(7,2.5){(b)}
\end{picture}
\end{figure}
By construction, the tableau $\psi(t)$ has shape $\theta/(k^{m-r})$ and content $(k^r)$.  It remains to show that $\psi(t)$ is indeed an LR-filling.  
Since the columns of $\psi(t)$ are strictly increasing, there can be
at most one $r-1$ in each column.  Since there are $k$ entries of $r-1$ in $\psi(t)$, every
$r$ in $\psi(t)$ must have a corresponding $r-1$ below it and weakly to the right.  
It follows that $\psi(t)$ is an LR-filling. 

Finally,  let $\theta=((2k)^r, k^{m-2r})$.  Since
  $$c^\theta_{(k^{m-r+1}), (k^{r-1})} = 0,\qquad {\rm while}\qquad c^\theta_{(k^{m-r}), (k^r)} = 1,$$
  we get an explicit verification 
that $d_{r-1} < d_r$.  
\end{proof}

\begin{remark}
We have just seen\footnote{The case $a=b$ has previously been obtained by Kirillov \cite{Kir84}. See Kleber \cite{Kle01} for a nice proof. } that, for all $k$, 
\[
     s_{k^a} s_{k^b} - s_{k^{a+1}} s_{k^{b-1}}
\]
is Schur-positive whenever $a \geq b$.
Now when $k=1$, this is nothing but the (dual) Jacobi-Trudi identity:
\[
s_{1^a} s_{1^b} - s_{1^{a+1}} s_{1^{b-1}} 
= \omega(s_{(a,b)}).
\]
For general $k$, we see that 
\begin{equation}\label{detspositive}
s_{k^a} s_{k^b} - s_{k^{a+1}} s_{k^{b-1}} = 
\omega \left( \left| \begin{array}{cc} s_{(a^k)} & s_{((a+1)^k)} \\ s_{((b-1)^k)} & s_{(b^k)} 
\end{array} \right| \right),
\end{equation}
so the determinant in \eqref{detspositive} is itself Schur-positive.
This observation has lead us to ask what happens if we take any Jacobi-Trudi determinant
   $$s_\mu = \det(s_{\mu_i -i+j})^p_{i,j=1} $$
and consider the corresponding ``exploded'' Jacobi-Trudi determinant
   $$\det(s_{((\mu_i -i+j)^k)})^p_{i,j=1}$$
for $k \geq 2$.
When is the result Schur-positive?  When $p=3$, we have systematically verified that these exploded determinants are Schur-positive whenever the resulting degree (as a formal power series) is at most 90.  
On the other hand, when $p=4$, there are simple examples of exploded
Jacobi-Trudi determinants that are not Schur-positive.  
\end{remark}

\section*{Acknowledgements}
The authors are grateful to Riccardo Biagioli and Sergey Fomin for helpful and interesting discussions.

\bibliography{../../master}
\bibliographystyle{plain}

\end{document}